\documentclass[11pt, dina4]{article}
\usepackage{graphicx,epsfig}
\usepackage{epstopdf}
\usepackage{amssymb,amsmath}
\usepackage{cases}
\usepackage{amsthm}
\usepackage{caption,lipsum}
\usepackage{stmaryrd}
\usepackage{tabularx}
\usepackage{subfigure}
\usepackage{color, soul}
\usepackage{multirow}
\usepackage{empheq,float} 
\DeclareGraphicsExtensions{.eps}
\usepackage{hyperref}
\usepackage[margin=1in]{geometry} 
\usepackage{comment}

\newcommand{\p}{\partial}
\newcommand{\Og}{\Omega}
\newcommand{\fl}[2]{\frac{#1}{#2}}

\newcommand{\nn}{\nonumber}

\newcommand{\Dt}{\Delta}
\newcommand{\ft}{\mathcal{F}}
\newcommand{\bxi}{{\boldsymbol \xi}}

\newcommand{\bea}{\begin{eqnarray}}
\newcommand{\eea}{\end{eqnarray}}
\newcommand{\beas}{\begin{eqnarray*}}
\newcommand{\eeas}{\end{eqnarray*}}
\newtheorem{remark}{Remark}[section]

\newcommand{\bx}{{\bf x} }
\newcommand{\by}{{\bf y} }
\newcommand{\md}{\mathrm{\,d}}

\newcommand{\bb}{\vskip 10pt}
\definecolor{ForestGreen}{rgb}{0.0, 0.5, 0.0}

\newcommand{\bT}{{\mathbb T}}

\newcommand{\fe}{\mathrm{e}}
\numberwithin{equation}{section}

\title{Fourier pseudospectral methods for the  variable-order space fractional wave equations}

\author{Yanzhi Zhang\thanks{Department of Mathematics and Statistics, Missouri University of Science and Technology, Rolla, MO 65409 (Email:  zhangyanz@mst.edu)}, \ \
Xiaofei Zhao\thanks{School of Mathematics and Statistics \& Computational Sciences Hubei Key Laboratory, Wuhan University, Wuhan, 430072, China (Email: matzhxf@whu.edu.cn)}, \ \
Shiping Zhou\thanks{Department of Mathematics and Statistics, Missouri University of Science and Technology, Rolla, MO 65409 (Email:  szb5g@mst.edu)}
}

\begin{document}
\date{}
\maketitle

\begin{abstract}
In this paper, we propose Fourier pseudospectral methods to solve the variable-order space fractional wave equation and develop an accelerated matrix-free approach for its effective implementation. 
In constant-order cases,  our methods can be efficiently implemented via the (inverse) fast Fourier transforms, and the computational cost at each time step is ${\mathcal O}(N\log N)$ with $N$  the total number of spatial points.  
However, this fast algorithm fails in the variable-order cases due to the spatial dependence of the Fourier multiplier. 
On the other hand, the direct matrix-vector multiplication approach becomes impractical due to excessive memory requirements. 
To address this challenge, we proposed an accelerated matrix-free approach for the efficient computation of variable-order cases. 
The computational cost is ${\mathcal O}(MN\log N)$ and storage cost ${\mathcal O}(MN)$, where $M \ll N$.  
Moreover, our method can be easily parallelized to further enhance its efficiency. 
Numerical studies show that our methods are effective in solving the variable-order space fractional wave equations, especially in high-dimensional cases.  
Wave propagation in heterogeneous media is studied in comparison to homogeneous counterparts. We find that wave dynamics in fractional cases become more intricate due to nonlocal interactions. Specifically, dynamics in heterogeneous media are more complex than those in homogeneous media. \\

\end{abstract}

\section{Introduction}
\setcounter{equation}{0}
\label{section1}

Recently, nonlocal fractional wave equations have garnered considerable attention in modeling wave propagation within complex media, including biomedical materials, irregular porous media, and fractal rock layers \cite{Autuori2018, Chen2004, Treeby2010, Zhu2014, Meerschaert2014, Maestas2016}. Specifically, variable-order space fractional wave equations have found extensive applications in studying wave dispersion and attenuation in heterogeneous media. 
In  \cite{Zhu2014, Chen2016, Yao2017, Mu2021} and many other studies, the seismic wave propagation in heterogeneous media is described by the following variable-order space fractional wave equations:
\begin{equation}\label{seismicwave}
\fl{1}{c^2(\bx)}\fl{\partial^2 u(\bx, t)}{\partial t^2}  = \eta(\bx) (-\Dt)^{1+\gamma(\bx)} u(\bx,t) + \tau(\bx)\fl{\partial}{\partial t}\big[(-\Dt)^{\fl{1}{2}+\gamma(\bx)}u(\bx,t)\big],
\end{equation}
where $u(\bx, t)$ is a wavefield function of space $\bx \in {\mathbb R}^d$ and time $t > 0$. 
For $s(\bx) > 0$, the operator $(-\Delta)^{s(\bx)}$ represents the variable-order fractional Laplacian,  which will be defined later. 
The  coefficient functions are given by \cite{Zhu2014, Chen2016, Yao2017, Mu2021}
\begin{align*}
c(\bx) &= c_0\cos\big(\pi\gamma(\bx)/2\big), \\   
\eta(\bx) &= -{c_0^{2\gamma(\bx)}}\omega_0^{-2\gamma(\bx)}\cos\big(\pi\gamma(\bx)\big), \\
\tau(\bx) &= -{c_0^{2\gamma(\bx)-1}}\omega_0^{-2\gamma(\bx)}\sin\big(\pi\gamma(\bx)\big),
\end{align*}
with $c_0$ representing the phase velocity at the reference frequency $\omega_0$.  
Here,  the function $\gamma(\bx) = {\rm arctan}(1/Q(\bx))/\pi$ with the quality factor $Q(\bx) > 0$.   
It is clear that $0 < \gamma(\bx) < 0.5$. 
If $\gamma = 0$, the model (\ref{seismicwave}) reduces to the classical acoustic wave equation. 
It shows in \cite{Meerschaert2014} that the variable-order fractional wave equations allow for varying attenuation indices in different regions to fully capture the anisotropic nature of wave propagation in complex media. 
Moreover, the viscoacoustic wave equation (\ref{seismicwave}) can effectively describe  the velocity dispersion and amplitude loss phenomena during the seismic waves  propagation. 
Consequently, it facilitates seismic imaging with reduced efforts in compensating for attenuation (see e.g., \cite{Zhu2014, Meerschaert2014, Chen2016, Yao2017, Mu2021} and other subsequent studies).

The variable-order fractional Laplacians find extensive applications in modeling heterogeneous properties of complex systems \cite{Chen2016, Zhu2014, Meerschaert2014, Sun2019, Wu2024, DElia2022, Cai2018, Yu2022}. 
It is defined as a pseudo-differential operator with symbol $|\bxi|^{2s(\bx)}$  \cite{Bass1988, Bass2004, Kikuchi1997, Leopold1999, Samko1993c, Sun2019, Jacob1993}:  
\begin{equation}\label{pseudo}
    (-\Dt)^{s(\bx)} u(\bx) = \int_{{\mathbb R}^d}  \widehat{u}(\bxi) |\boldsymbol\xi|^{2s(\bx)} \fe^{2\pi i \bx\cdot{\boldsymbol\xi}}\md\bxi, \qquad \text{ for} \  s(\bx) > 0,
\end{equation}
 where $\widehat{u}(\bxi)$ represents the Fourier transform of $u(\bx)$.  
It is assumed assume that $s(\bx) \ge \inf s(\bx) > 0$, and $s(\bx)$ is H\"older continuous. 
If $s(\bx) \equiv s $  is constant, the operator in (\ref{pseudo}) collapses to the celebrated (constant-order) fractional Laplacian, which can be simplified to
\begin{equation}\label{pseudo1}
    (-\Dt)^{s} u(\bx) = \ft^{-1}\big[|\bxi|^{2s} \ft[u]\big], \qquad \text{ for} \ s > 0,
\end{equation}
where $\ft$ and $\ft^{-1}$ represent the Fourier transform and its associated inverse transform, respectively.  
If the exponent $0 < s(\bx) < 1$, the variable-order fractional Laplacian can be also defined in a hypersingular integral form \cite{Bass1988, Chen2020, Kuhn2021, Samko1993c}: 
\bea\label{integralFL}
(-\Dt)^{s(\bx)} u(\bx) = \fl{4^{s(\bx)}s(\bx)\Gamma\big(s(\bx) + \fl{d}{2}\big)}{\sqrt{\pi^d}\,\Gamma\big(1-s(\bx)\big)}\, {\rm P. V.}\int_{{\mathbb R}^d}\fl{u(\bx) - u(\by)}{|\bx - \by|^{d+2s(\bx)}}\, d\by,  
\eea
where P.V. stands for the principal value integral, and  $\Gamma(\cdot)$ is the Gamma function.  
More discussion of the variable-order fractional Laplacian $(-\Dt)^{s(\bx)}$ can be found in \cite{Samko1993, Wu2024} and references therein. 

Currently, numerical studies on the variable-order space fractional wave equation still remain limited.  
The main numerical challenge lies in the lack of effective numerical methods for computing the variable-order fractional Laplacian $(-\Dt)^{s(\bx)}$. 
Compared to its constant-order counterpart, the combination of nonlocality and heterogeneity in variable-order fractional Laplacian introduces significant storage and computational challenges.  Consequently, many numerical methods developed for the constant-order fractional Laplacian become ineffective for computing the variable-order cases. 
So far, two methods have been recently proposed to compute the variable-order fractional Laplacian: one is meshfree radial basis function methods \cite{Wu2024}, and the other is finite element methods \cite{DElia2022}. 
The lack of numerical methods for variable-order fractional Laplacian severely impedes numerical studies of the fractional wave equations (\ref{model}), particularly in higher dimensions ($d > 1$).  
In \cite{Zhu2014},  the fractional wave equation is studied by approximating the variable-order fractional Laplacian with an averaged constant order $\bar{s} = {\rm avg}\{s(\bx)\}$. 
Later, a weighted sum of multiple constant-order fractional Laplacians is introduced in \cite{Li2015, Mu2021} to approximate the variable-order fractional Laplacian. 
These strategies that approximate the variable-order fractional Laplacian with its constant-order counterpart might alleviate numerical challenges in solving the fractional wave equation. 
However, they may also introduce potential challenges in accurately describing the heterogeneity.

In this paper, we introduce Fourier pseudospectral methods to solve the variable-order space fractional wave equation (\ref{model}) and develop an accelerated matrix-free approach for its effective implementation. 
To  facilitate our discussion, we will focus on the variable-order space fractional wave equation of the general form \cite{Zhu2014, Meerschaert2014, Sun2019, Yao2017}: 
\begin{equation}\label{model}
\begin{aligned}
\partial_{tt}u(\bx,t)= -\kappa(-\Dt)^{s(\bx)}u(\bx,t) + f(u), \qquad \ &  \mbox{for} \  \bx\in\bT^d,\quad t>0,\\
u(\bx,0)=\phi(\bx),\quad \partial_tu(\bx,0)=\psi(\bx), \qquad \ &\mbox{for} \ \bx\in\bT^d,
\end{aligned}
\end{equation}
where constant $\kappa > 0$, and $\bT^d$ denotes a $d$-dimensional torus for $d=1, 2$, or $3$. 
The model (\ref{model}) covers a broad class of wave equations. 
For $s(\bx)\equiv 1$, it reduces to the classical wave equations, while it becomes the biharmonic wave equation for $s(\bx)\equiv 2$ \cite{Roetman1967, Li2022}. 
In this study, we will focus on the power $0 < s(\bx) < 2$. 

Our Fourier pseudospectral methods uniformly solve both constant-order and variable-order problems. Moreover, if $s(\bx) \equiv s$ is a constant,  the (inverse) fast Fourier transform (FFT) can be utilized for its efficient implementation at a computational cost of ${\mathcal O}(N\log N)$ with $N$ the total number of spatial points. 
However, in variable-order cases, the spatial dependency of $s(\bx)$ causes the failure of inverse FFTs, rendering the fast algorithms designed for constant-order cases ineffective. 
On the other hand, the direct implementation of Fourier pseudospectral methods requires storing a full matrix and computing matrix-vector products at each time step.   
This leads to prohibitive storage and computational costs, particularly in high-dimensional cases, making it impractical. 
Hence, we propose an accelerated matrix-free approach for the efficient computation with computational cost ${\mathcal O}(MN\log N)$ and storage cost ${\mathcal O}(MN)$, where $M \ll N$.  
Numerical studies show that our accelerated matrix-free approach significantly outperforms the direct approach.  
For temporal discretization, we introduce and compare three methods, including Crank-Nicolson, leap-frog, and time splitting methods.
All of them have the second-order accuracy. 
The leap-frog and time splitting methods are explicit and thus are easier and more cost-effective to implement. 
The Crank--Nicolson method requires more computing time at each time step, but it is more stable. 
Compared to the other two methods, increasing the accuracy of the time-splitting method is much easier.  
Numerical experiments are performed to compare their accuracy and computing time in solving the fractional wave equations. 
The interactions of solitary waves and the dispersion of waves are studied in both homogeneous (constant $s$) and heterogeneous (variable $s(\bx)$) media to understand the heterogeneity effects. 
Finally, we also apply our method to study the dispersion and attenuation of seismic waves in \cite{Zhu2014, Yao2017}. 

The paper is organized as follows. 
In Section \ref{section2}, we introduce a Fourier pseudospectral method for spatial discretization and propose an accelerated matrix-free approach to tackle the computational challenges caused by spatial heterogeneity. 
In Section \ref{section3}, we present and compare three temporal discretization schemes. 
In Section \ref{section4}, we conduct numerical experiments to examine the performance of our methods and explore wave dynamics in both homogeneous (constant $s$) and heterogeneous (variable $s(\bx)$) media. 
Conclusions and discussion are presented in Section \ref{section5}.

\section{Fourier pseudospectral methods}
\label{section2}

In this section, we focus on the spatial discretization of the variable-order fractional wave equation, and its temporal discretization will be discussed in Section \ref{section3}. 
The definition in (\ref{pseudo}) provides a unified pseudo-differential representation of the Laplace operator $(-\Dt)^{s(\bx)}$ for any $s(\bx) > 0$. 
Noticing this fact,  we apply the Fourier pseudospectral method for spatial discretization and develop a unified approximation for both constant-order and variable-order Laplacians. 
On the other hand, even though the constant- and variable-order Laplacians share the same discretization, numerical implementation of the variable-order cases is significantly more challenging (e.g. see comparison in \eqref{mu-eq}--\eqref{mu-eq1}). 
If $s(\bx)$ is spatially varying, the heterogeneity of operator $(-\Dt)^{s(\bx)}$ demands considerably more storage and computational costs and makes its numerical evaluation a formidable challenge, especially in high dimensions.

In the following, we first introduce our spatial discretization method, then compare numerical implementation of the constant-order and variable-order cases, and propose a matrix-free acceleration approach for efficient computation in the variable-order cases.  
For notational simplicity, let's first focus on the one-dimensional ($d = 1$) case with domain $\Og = (a, b)$ and introduce our method for the  problem: 
\bea\label{fw1d}
\begin{aligned}
    \p_{tt}u(x,t) = -\kappa (-\p_{xx})^{s(x)}u(x,t) + f(u), &\qquad \mbox{for}\ \  a \le x \le b, \  \quad  t>0,\\
u(x,0)=\phi(x), \quad \partial_t u(x,0)=\psi(x), &\qquad \mbox{for}\ \  a \le x \le b,
\end{aligned}
\eea
with periodic boundary conditions. 
The generalization of our method to higher dimensions (i.e., $d > 1$) will be discussed in Section \ref{section2-1}.  
Choose an even integer $J > 0$, and let mesh size $h =(b-a)/J$. 
Define the spatial grid points $x_j = a+jh$,  for $j = 0, 1, \ldots, J$. 
Assume the solution ansatz of (\ref{fw1d}) takes the form: 
\begin{equation}\label{ansatz}
        u^h(x, t) = \sum_{k=-J/2}^{J/2-1}\widehat{u}_k(t)\,\fe^{i\mu_k(x - a)},
\end{equation}
where we denote
\beas
 \mu_k=\frac{2\pi k}{b-a}, \quad \   \widehat{u}_k(t)=\frac{1}{J} \sum_{j = 0}^{J-1}u_j(t)\,\fe^{-i\mu_k(x_j-a)}, \qquad \mbox{for} \ \ -\fl{J}{2} \le k \le \fl{J}{2}-1,
\eeas
with $u_j(t)$ representing the numerical solution to $u(x_j, t)$.  
Substituting (\ref{ansatz}) into (\ref{pseudo}), we then obtain the numerical approximation of $(-\p_{xx})^{s(x)}u(x, t)$ at point $x = x_j$: 
\begin{equation}\label{FL1Ddis}
(-\p_{xx})^{s(x_j)}_h u_j(t) := \Big((-\p_{xx})_h^{s(x)}u(x, t)\Big)\Big|_{x = x_j}= \sum_{k=-J/2}^{J/2-1} |\mu_k|^{2s(x_j)}\,\widehat{u}_k(t)\,\fe^{i\mu_k(x_j-a)},
\end{equation}
for \,$j = 0, 1,\,\ldots,\,J-1$.  

Note that the Fourier pseudospectral approximation in (\ref{FL1Ddis}) holds for both constant-order and variable-order Laplacians. 
However, numerical implementation for spatially dependent $s(x)$ is considerably more challenging than that with a constant $s$. 
To see these challenges, let's assume that coefficients $\widehat{u}_k(t)$ in (\ref{FL1Ddis}) are known for $-J/2 \le k \le J/2-1$.  
If $s(x) \equiv s$ is a constant,  the approximation of the constant-order Laplacian $(-\p_{xx})^s$ in (\ref{FL1Ddis}) can be efficiently computed via the inverse FFTs with a computational cost of ${\mathcal O}(J\log J)$. 
In contrast,  if $s(x)$ is spatially varying,  the variable-order fractional Laplacian denotes a space-frequency mixed operator, and the spatial dependence of $|\mu_k|^{2s(x)}$ makes the inverse FFTs fail to calculate the summation in (\ref{FL1Ddis}). 
This difference could potentially affect the temporal approximation of the wave equation. 
To further show it, let's use the linear wave equation (i.e.  (\ref{fw1d}) with $f(u) = 0$) as an example. 
Substitute the ansatz (\ref{ansatz}) into the linear wave equation,  and then apply  FFT at both sides of it. 
\begin{itemize}
\item In \textit{homogenous} (i.e. constant $s$) cases,  we obtain a system of decoupled ordinary differential equations (ODEs) for $\widehat{u}_k(t)$, i.e.,
\begin{equation}\label{mu-eq}
 \fl{\md^2 \widehat{u}_k(t)}{\md t^2}+\kappa |\mu_k|^{2s}\widehat{u}_k(t)=0, \quad \  \mbox{for} \ \,  -\fl{J}{2} \le k \le \fl{J}{2}-1,
\end{equation}
with $s > 0$ a constant. 
The ODEs in (\ref{mu-eq}) can be solved independently. 
Particularly, they can be integrated in time exactly.
\item In contrast, the situation of \textit{heterogeneous} (i.e. variable $s(x)$) cases is more complicated, where we obtain a system of  ODEs for $\widehat{u}_k(t)$ as: 
\bea\label{mu-eq1}
 \fl{\md^2 \widehat{u}_k(t)}{\md t^2} =-\frac{\kappa}{J}\sum_{l=-J/2}^{J/2-1}\widehat{u}_l(t)\bigg(\sum_{j=0}^{J-1}|\mu_l|^{2s(x_j)}
\,\fe^{i(\mu_l-\mu_k)(x_j-a)}\bigg), \quad \mbox{for} \ \,  -\fl{J}{2} \le k \le \fl{J}{2}-1.
\eea
It shows that if $s(x)$ is spatially dependent,  the fast Fourier transform of the linear wave equation does not produce decoupled ODEs as in (\ref{mu-eq}),  but instead a coupled system for $\widehat{u}_k(t)$. 
Clearly, it is more challenging to solve the system in (\ref{mu-eq1}). 
\end{itemize}

\subsection{Direct matrix-vector approach} 
\label{section2-1}

The above  comparison and discussion suggest that the essential challenge in calculating the variable-order fractional Laplacian (\ref{FL1Ddis}) stems from the space-frequency mixed symbol. 
If $s(x)$ is spatially dependent, it is not beneficial to approximate the variable-order fractional Laplacian from $\widehat{u}_k(t)$ due to the failure of using inverse FFTs. 
Hence, one {\it direct matrix-vector approach} is to approximate the fractional Laplacian in space domain.  
By combining (\ref{ansatz}) and (\ref{FL1Ddis}), we can rewrite the approximation in (\ref{FL1Ddis}) as: 
\bea\label{FL1Ddis0}
(-\p_{xx})_{h}^{s(x_j)}u_j(t) = \frac{1}{J}\sum_{l = 0}^{J-1}\sum_{k=-J/2}^{J/2-1}|\mu_k|^{2s(x_j)}\,u_l(t)\,\fe^{i\mu_k(x_j-x_l)},\quad \mbox{for} \  \, 0 \le j \le J-1.  
\eea
Compared to (\ref{FL1Ddis}), the formulation in (\ref{FL1Ddis0}) enables us to approximate the fractional Laplacian directly from function $u_j(t)$, and thus avoids the computational cost in obtaining $\widehat{u}_k(t)$.  

Denote vector ${\bf u}(t):=\big(u_0(t), \, u_1(t),\,\cdots,\, u_{J-1}(t)\big)^{T}$.  
The approximation  in (\ref{FL1Ddis0}) can be formulated into a matrix-vector form, i.e. 
$(-\p_{xx})_{h}^{s(x)} u(x, t) = A_{s(x)}{\bf u}(t)$, where the entries of  matrix  $A_{s(x)} = \big\{a_{jl}\big\}_{J\times J}$ are given by
\bea\label{mat-rep}
&&a_{jl}=\frac{1}{J}\sum_{k=-J/2}^{J/2-1}|\mu_k|^{2s(x_j)}\,\fe^{i\mu_k(x_j-x_l)} \nn\\
&&\hspace{0.6cm}= \fl{2}{J} \bigg((-1)^{j-k} + \sum_{k = 1}^{J/2-1} |\mu_k|^{2s(x_j)}\cos(\mu_k(x_j-x_l))\bigg), \quad \  \mbox{for} \ \, 0 \le j,\,l \le J-1.\qquad \ 
\eea
Generally,  $A_{s(x)}$ is a full matrix, which requires a memory cost of ${\mathcal O}(J^2)$. 
Due to the spatially varying $s(x)$, the computational cost in assembling the matrix is ${\mathcal O}(J^3)$, while  the cost in calculating matrix-vector product is ${\mathcal O}(J^2)$.  
Furthermore, the condition number of $A_{s(x)}$ increases dramatically as the number of points $J$ increases. 
In special cases of constant order $s(x) \equiv s$,  $A_{s(x)}$ reduces to a symmetric Toeplitz matrix. 
Consequently, the memory cost reduces to ${\mathcal O}(J)$.  
Moreover, its Toeplitz structure enables us to design fast algorithms for matrix-vector multiplication via fast Fourier transforms (FFTs) at a cost of ${\mathcal O}(2J\log(2J))$ \cite{Duo2019-FDM}. 

Due to the high memory and computational costs, generalizing the direct matrix-vector method to higher dimensions becomes challenging and impractical, particularly as the number of points is large; 
see more discussion and illustration in Section \ref{section4-2}.

\subsection{Accelerated matrix-free approach} 
\label{section2-2}

Next,  we propose an {\it accelerated matrix-free approach} to essentially minimize the storage and computational costs in computing the variable-order fractional Laplacian from (\ref{FL1Ddis}) and (\ref{FL1Ddis0}). 
As discussed earlier, the failure of using inverse FFTs to accelerate the calculation of (\ref{FL1Ddis}) is caused by the spatial dependence of the Fourier multiplier $|\mu_k|^{2s(x)}$. 
Hence to resolve it, we focus on separating the spatial dependence from the Fourier multiplier. 
Consider the exponential function $a^z$ for $z \in {\mathbb R}$ and the base $a > 0$. 
Its Taylor expansion at point $z = z_0$ gives
\bea\label{a-expansion}
a^z = a^{z_0} \sum_{m = 0}^\infty\,\fl{1}{m!}\,(z - z_0)^m\,(\ln a)^m. 
\eea
Similarly,  we can apply the Taylor expansion to the Fourier multiplier $|\mu_k|^{2s(x)}$.   
Letting $a = |\mu_k|^{2}$,\, $z = s(x)$,\, and $z_0 = s_0$ in (\ref{a-expansion}), we obtain 
\bea\label{mu-expansion}
|\mu_k|^{2s(x)} =  |\mu_k|^{2s_0} + |\mu_k|^{2s_0}\bigg(\sum_{m=1}^{\infty}\frac{1}{m!}\big[s(x)-s_0\big]^m \big(\ln{|\mu_k|^{2}}\big)^m\bigg),\qquad \mbox{for} \ \ k\neq0.
\eea 
Note that for $k = 0$, no expansion is needed as $\mu_0 = 0$.   
The expansion in (\ref{mu-expansion}) implies that the variable-order fractional Laplacian $(-\Dt)^{s(x)}$ can be viewed as a perturbed form of the constant-order fractional Laplacian $(-\Dt)^{s_0}$, and the  perturbation is spatially dependent.   
Here, we choose the constant $s_0 = \big[\int_a^b s(x)  \md x\big]/(b-a)$.  

Taking a sufficiently large $M$, we can truncate the summation in (\ref{mu-expansion}) into $M$ terms.  
Substituting the truncated $M$-term expansion of (\ref{mu-expansion}) into (\ref{FL1Ddis}) yields our new approximation to the variable-order fractional Laplacian:
\bea\label{Lu}
(-\Dt)_{h}^{s(x_j)}u_j(t) = \sum_{m=0}^{M}\big[s(x_j)- s_0\big]^m
\Bigg(\sum_{\substack{k=-J/2 \\ k \neq 0}}^{J/2-1}\frac{\big(\ln{|\mu_k|^2}\big)^m}{m!}\big[|\mu_k|^{2s_0}\widehat{u}_k(t)\big]\fe^{i\mu_k(x_j-a)}\Bigg), 
\eea
for $j = 0, 1, \, \ldots, \, J-1$.  
Note that the term with $m = 0$ provides an approximation of the constant-order fractional Laplacian $(-\Delta)^{s_0}$. 
Our scheme in (\ref{Lu}) provides a matrix-free approximation to the variable-order fractional Laplacian. 
The inner summation over $k$  can be efficiently calculated by inverse FFTs with a computational cost of ${\mathcal O}(J \log J)$. 
Consequently, the total computational cost of (\ref{Lu}) in approximating the variable-order fractional Laplacian is ${\mathcal O}(M J \log J)$. 
Usually the number $M \ll J$, and thus the computational cost in calculating (\ref{Lu}) is significantly lower than ${\mathcal O}(J^2)$ of the direct matrix-vector approach in (\ref{FL1Ddis0}). 
Moreover, the scheme  (\ref{Lu}) does not require assembling and storing matrices.  

\begin{figure}[htb!]
\centerline{
(a)\includegraphics[width=6.86cm, height=5.06cm]{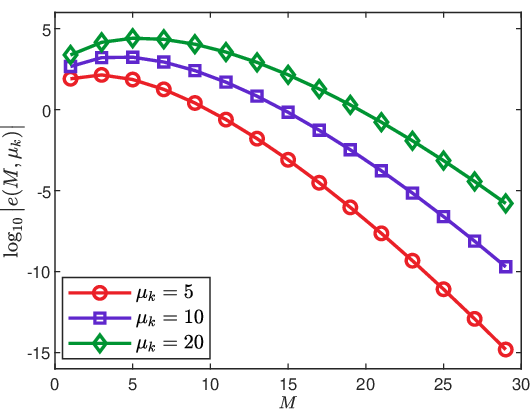}
\hspace{2mm}
(b)\includegraphics[width=6.86cm, height=5.06cm]{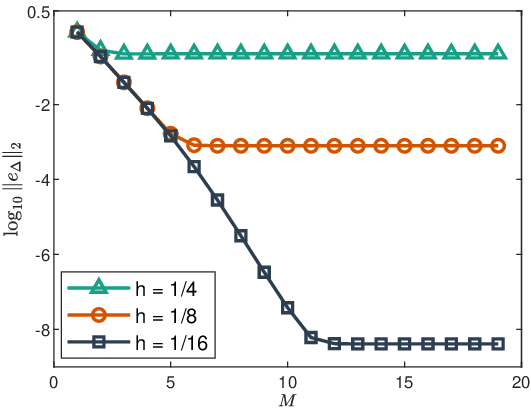}}
\caption{(a) Illustration of the relation between $|e(M, \mu_k)|$, $M$ and $\mu_k$, where $s_0 = 1$ is used. (b) The $l_2$-norm errors in approximating function $(-\Dt)^{s(x)}u(x)$ for different $h$ and $M$, where $s(x) = 1+0.3\sin(\pi x/8)$ and $u(x)$ is defined in (\ref{example3}).} \label{Figure0}
\end{figure}

Generally, the choice of $M$ in scheme (\ref{Lu}) depends on exponent $s(x)$, mesh size $h$,  and the solution behavior in frequency domain.  
Denote
\beas
e(M, \mu_k) = \frac{2^M\big(\ln{|\mu_k|}\big)^{M}}{M!}\,|\mu_k|^{2s_0}, \qquad\mbox{for $M \in {\mathbb N}$. } 
\eeas
It is easy to see that $e(M, \mu_k)$ decreases quickly with $M$ increasing; see illustration in Figure \ref{Figure0} (a). 
The larger the number $M$, the smaller the truncation errors of (\ref{Lu}), as $|s(x) - s_0| < 1$.  
Moreover, notice that  $e(M, \mu_k)$ serves as the coefficient of $\widehat{u}_k(t)$, while solution $\widehat{u}_k(t)$ usually decreases quickly with respect to $|\mu_k|$ and has a finite bandwidth centered at $\mu_k = 0$. 
This implies that even though $e(M, \mu_k)$ increases with $\mu_k$, the impact of truncation becomes negligible if $|\mu_k|$ is outside of the bandwidth. 
Figure \ref{Figure0} (b) further shows the numerical errors in approximating the variable-order fractional Laplacian for different $M$.  
For a fixed mesh size $h$, numerical errors first decrease as $M$ increases.
However, when $M$ is sufficiently large, the discretization errors become dominant, and the truncation effect can be ignored. 
Furthermore, as the mesh size $h$ decreases, the required terms (i.e. $M$) increase. 
More discussion and comparison of different $M$ can be found in Section \ref{section4-2}. 

\subsection{Generalization to higher dimensions}
\label{section2-1}

Our matrix-free scheme in (\ref{Lu}) provides an effective approach to compute the variable-order fractional Laplacian. 
Moreover, it can be easily generalized into higher dimensions.  
For the convenience of readers, we present the generalized scheme for $d \ge 1$ as follows. 

Let the $d$-dimensional domain $\Og = (a_1, b_1)\times(a_2, b_2)\times \cdots \times (a_d, b_d)$.  
For $1 \le m \le d$, choose even integers $J_m > 0$,  and define mesh size  $h_m = ({b_m-a_m})/{J_m}$. 
For notational simplicity,  we denote the index sets 
\beas
&&S_x = \{{\bf j} = (j_1,\, j_2, \,\cdots,\, j_d), \ \ \mbox{for} \ \,  0 \le j_m \le J_m-1 \}, \\
&&S_\xi = \{{\bf k} = (k_1,\, k_2, \,\cdots,\, k_d), \ \ \mbox{for} \ -J_m/2 \le k_m \le J_m/2-1, \ \mbox{but} \ \, {\bf k} \neq {\bf 0} \},\qquad\qquad
\eeas
and $S_\xi^0 = S_\xi \cup \{ {\bf k} = {\bf 0} \}$. 
The $d$-dimensional  grid points are denoted as $\bx_{\bf j} = \big(x_{j_1}^{(1)}, \, x_{j_2}^{(2)}, \, \cdots, x_{j_d}^{(d)}\big)$ for \,${\bf j} \in S_x$,  and ${\boldsymbol\mu_{\bf k}} = \big(\mu_{k_1}^{(1)},\,\mu_{k_2}^{(2)},\, \cdots, \mu_{k_d}^{(d)}\big)$ for ${\bf k} \in S_\xi$, where $x_l^{(m)} = a_m + lh_m$ and  $\mu_l^{(m)}= 2\pi l/(b_m - a_m)$.  
Let $u_{\bf j}(t)$ represent the numerical solution of $u(\bx_{\bf j}, t)$, and assume
\begin{equation}\label{ansatz-d}
u_{\bf j}(t) =  \sum_{{\bf k}\, \in \, {S}_\xi^0}\widehat{u}_{\bf k}(t)\, \fe^{i{\boldsymbol\mu}_{\bf k}{\boldsymbol\cdot}(\bx_{\bf j} - {\bf a})},  \qquad\mbox{for} \ \,{\bf j} \in S_x,
\end{equation}
where ${\boldsymbol\mu}_{\bf k}{\boldsymbol\cdot}(\bx_{\bf j} - {\bf a}) = \prod_{m = 1}^d \mu_{k_m}^{(m)} \big(x_{j_m}^{(m)} - a_m\big)$, and 
\beas
\widehat{u}_{\bf k}(t)=\Big(\prod_{m=1}^d\frac{1}{J_m}\Big)\sum_{{\bf j}\,\in\,{S}_x}u_{\bf j}(t)\,\fe^{-i{\boldsymbol\mu}_{\bf k}{\boldsymbol\cdot}(\bx_{\bf j} - {\bf a})},   \qquad\mbox{for} \ \,{\bf k} \in S_\xi^0. 
\eeas
Combining the solution ansatz in (\ref{ansatz-d}) with the definition in (\ref{pseudo}) leads to the numerical approximation of $d$-dimensional fractional Laplacian $(-\Dt)^{s(\bx)}$ as
\bea\label{FLdis-d}
(-\Dt)^{s(\bx_{\bf j})}_h u_{\bf j}(t) &:=& \Big((-\Dt)_h^{s(\bx)}u(\bx, t)\Big)\Big|_{\bx = \bx_{\bf j}} \nn\\
&=& \sum_{{\bf k} \in\,{S}_\xi} |\boldsymbol\mu_{\bf k}|^{2s(\bx_{\bf j})}\,\widehat{u}_{\bf k}(t)\,\fe^{i {\boldsymbol \mu}_{\bf k}{\boldsymbol\cdot}(\bx_{\bf j} - {\bf a})}, \qquad \mbox{for} \ \ {\bf j} \in S_x,\qquad 
\eea
where $|\boldsymbol\mu_{\bf k}|^{2s(\bx_{\bf j})} = \big[\sum_{m=1}^d \big|\mu_{k_m}^{(m)}\big|^2\big]^{s(\bx_{\bf j})}$. 
Similar to the one-dimensional cases, if $s(\bx) \equiv s$ is a constant, the approximation in (\ref{FLdis-d}) can be efficiently computed using the $d$-dimensional inverse FFT. 

However,  when $s(\bx)$ varies spatially, the application of inverse FFTs fails to compute the summation in (\ref{FLdis-d}). 
While directly calculating this summation leads to a computational cost of  ${\mathcal O}(N^2)$ with $N = \prod_{m=1}^d J_m$ the total number of spatial points.  
To reduce the computational complexity, we adopt the similar idea as in (\ref{mu-expansion}) to expand the Fourier multiplier and thus obtain the further approximation: 
\bea\label{Lu-d}
(-\Dt)_{h}^{s(\bx_{\bf j})}u_{\bf j}(t) = \sum_{m=0}^{M}\big[s({\bf x_j})- s_0\big]^m\Bigg(\sum_{{\bf k}\in{S_\xi}}\frac{\big(\ln{|\boldsymbol\mu_{\bf {\bf k}}|^2}\big)^m}{m!}\big[|\boldsymbol\mu_{\bf k}|^{2s_0}\widehat{u}_{\bf k}(t)\big]\fe^{i{\boldsymbol\mu}_{\bf k}{\boldsymbol\cdot}(\bx_{\bf j}-{\bf a})}\Bigg),
\eea
where $s_0$ is chosen as  $s_0 = \big(\int_\Og s(\bx) \md\bx\big)/|\Og|$. 
The coefficients of $\widehat{u}_{\bf k}(t)$ in the inner summation of (\ref{Lu-d}) are independent of space $\bx$. 
Hence, the inner summation can be efficiently computed by inverse  FFTs with computational cost of ${\mathcal O}(N \log N)$. 
As a result, the computational cost for approximating the variable-order fractional Laplacian in (\ref{Lu-d})  is ${\mathcal O}(M N \log N)$ with $M \ll N$.  
The scheme (\ref{Lu-d}) provides an effective approach to compute the $d$-dimensional fractional Laplacian $(-\Dt)^{s(\bx)}$. 
Note that the direct approach (e.g., (\ref{FL1Ddis0}) for $d = 1$) becomes impractical in higher dimensions due to the formidable challenges in storing the matrix and performing matrix-vector multiplications. 

In summary, if $s(\bx) \equiv s$ is a constant, the fractional Laplacian can be efficiently computed either in the frequency domain from $\widehat{\bf u}(t)$ or in the space domain from ${\bf u}(t)$.   
Let's take the one-dimensional case as an example. 
The computation in the frequency domain via (\ref{FL1Ddis}) can be directly realized by inverse FFTs at a cost of ${\mathcal O}(N\log N)$.  
\begin{table}[htb!]
\begin{center}
\begin{tabular}{|l|cc|cc|}
\hline 
&  Constant $s$ & & & Variable $s(\bx)$ \\
\hline
Direct method in space domain &  ${\mathcal O}(2N\log(2N))$ & & &$\mathcal{O}(N^2)$ \\
Scheme (\ref{FLdis-d}) in frequency domain & ${\mathcal O}\big(N\log N)$  & & & $\mathcal{O}(N^2)$ \\
Scheme (\ref{Lu-d}) in frequency domain  & n.a. & \quad &  & $\mathcal{O}(MN\log N)$ \\
\hline
\end{tabular}
\caption{Summary of the numerical scheme and corresponding computational costs in approximating the constant-order and variable-order fractional Laplacians. Note that the direct approach (\ref{FL1Ddis0}) in space domain can be straightforwardly generalized for $d > 1$. }\label{Tab00}
\end{center}
\end{table}
While the implementation in the space domain via (\ref{FL1Ddis0})  involves matrix-vector product. 
In this case, the matrix is symmetric  Toeplitz, and fast algorithms can be designed for matrix-vector multiplication using FFTs with computational cost of  ${\mathcal O}(2N\log(2N))$.  
Generally, the computation of the variable-order fractional Laplacian is more complicated and costly. 
For the convenience of readers, we summarize and compare different approaches in Table \ref{Tab00}.  
In the variable-order cases, the storage cost is ${\mathcal O}(N^2)$ for both the direct method in space domain and the scheme (\ref{FLdis-d}) in frequency domain. 

\section{Time discretization}
\label{section3}

In the previous section, we introduce a Fourier pseudospectral method for spatial discretization and propose an accelerated matrix-free approach (\ref{Lu-d})  to efficiently compute the variable-order fractional Laplacian. 
In this section, we focus on the temporal discretization and discuss three numerical methods.   
Denote $\tau > 0$ as the time step, and define time sequence $t_n = n \tau$ for $n = 0, 1, \ldots$. 
Let $u_{\bf j}^n$ represent the numerical approximation to solution $u(\bx_{\bf j}, t_n)$, for ${\bf j} \in S_x$ and $n = 0, 1, \ldots$. 

\subsection{Crank--Nicolson Fourier pseudospectral method }
\label{section3-1}

Let's start with the semi-discretization of the wave equation (\ref{model}), i.e., 
\bea\label{semi}
\fl{\md^2 u_{\bf j}(t)}{\md t^2} = -\kappa (-\Dt)_h^{s(\bx_{\bf j})} u_{\bf j}(t) + f(u_{\bf j}(t)), \quad  \ \mbox{for} \ \ {\bf j} \in S_x,
\eea
where $(-\Dt)_h^{s(\bx_{\bf j})}$ represents the numerical approximation of the Laplace operator at point $\bx = \bx_{ j}$. 
The application of Crank--Nicolson method to (\ref{semi}) yields the fully discretized scheme as 

\bea\label{CNscheme}
\frac{u_{\bf j}^{n+1}-2u_{\bf j}^n+u_{\bf j}^{n-1}}{\tau^2}=-\fl{\kappa}{2} (-\Dt)_h^{s(\bx_{\bf j})}  \Big(u_{\bf j}^{n+1} + u_{\bf j}^{n-1}\Big) + \fl{f(u_{\bf j}^{n+1}) + f(u_{\bf j}^{n-1})}{2},
\eea
for ${\bf j} \in S_x$ and $n = 1, 2, \ldots$. 
For $n = 0$, we can get $u_{\bf j}^0$ exactly from the initial condition $u(\bx, 0)$. 
To obtain $u_{\bf j}^1$, we first take the Taylor expansion of $u(\bx, t_1)$ at $t = 0$, i.e., 
\beas
u(\bx, t_1) = u(\bx, 0) + \tau \p_t u(\bx, 0) + \fl{\tau^2}{2}\p_{tt} u(\bx, 0) + O(\tau^3), 
\eeas
and then substitute the initial conditions and the wave equation (\ref{model}) at $\bx = \bx_{\bf j}$ into it. 
Hence, we obtain the approximation at $n = 0$ and $1$ as: 
 \begin{subequations}\label{IC-dis}
 \begin{align}
&u_{\bf j}^0 = \phi(\bx_{\bf j}), \\
&u_{\bf j}^1=\phi(\bx_{\bf j})+\tau\,\psi(\bx_{\bf j}) + \frac{\tau^2}{2} \bigg(-\kappa (-\Dt)_h^{s(\bx_{\bf j})}\phi(\bx_{\bf j}) + f\big(\phi(\bx_{\bf j})\big)\bigg), \quad \ \mbox{for} \ \ {\bf j}\in S_x. 
\end{align}
\end{subequations}
We will refer to scheme  (\ref{CNscheme})--(\ref{IC-dis})  as the Crank--Nicolson Fourier pseudospectral  (CNFP)  method. 
It has the second-order temporal accuracy and spectral spatial accuracy. 
Moreover, the CNFP scheme is unconditionally stable. 
It is implicit in time, and at each time step $t = t_n$ the resulting nonlinear system is solved by the conjugate gradient method.

Tables \ref{table-CN} shows the temporal errors and convergence rate of the CNFP method in solving the one-dimensional nonlinear fractional wave equation in (\ref{example1}), while the spatial errors will be presented in Section \ref{section4}. 
The numerical parameters are the same as those used in Example 1 of Section \ref{section4}. 
Since the exact solution is unknown, we use numerical solution with fine mesh $h = 1/64$ and time step $\tau = 0.0001$ as the reference in computing numerical errors. 
\begin{table}[htb!]
\begin{center}
\begin{tabular}{|c|c|c|c|c|c|c|c|c|c|c|}
\hline
\multirow{2}{*}{$\tau$} &\multicolumn{2}{|c|}{$s(x) \equiv 0.5$} & \multicolumn{2}{|c|}{$s(x) \equiv 1$} &  \multicolumn{2}{|c|}{$s(x) \equiv 1.3$} &  \multicolumn{2}{|c|}{$s_1(x)$} &  \multicolumn{2}{|c|}{$s_2(x)$} \\
\cline{2-11}
& error & c.r &  error & c.r. & error & c.r. & error & c.r. & error & c.r. \\
\hline
$2^{-7}$  & $5.618e$-6 & --   &  2.515$e$-5 & --  & 6.978$e$-5 & -- & 2.803$e$-5 & -- & 3.630$e$-5 & -- \\ 
$2^{-8}$  & $1.408e$-6 & 2.00 &  6.297$e$-6 & 2.00 & 1.748$e$-5 & 2.00 & 7.021$e$-6 & 2.00 & 9.092$e$-6 & 2.00 \\
$2^{-9}$  & 3.500$e$-7 & 2.01 &  1.559$e$-6 & 2.01 & 4.329$e$-6 & 2.01 & 1.739$e$-6 & 2.01 & 2.251$e$-6 & 2.01 \\ 
$2^{-10}$ & 8.490$e$-8 & 2.05 &  3.714$e$-7 & 2.07 & 1.032$e$-6 & 2.07 & 4.143$e$-7 & 2.07 & 5.366$e$-7 & 2.07 \\ 
$2^{-11}$ & 1.941$e$-8 & 2.13 &  7.410$e$-8 & 2.32 & 2.065$e$-7 & 2.32 & 8.281$e$-8 & 2.32 & 1.073$e$-7 & 2.32 \\ \hline
\end{tabular}
\caption{Temporal errors { $\|u(t) - u^{h, \tau}(t)\|_{l^2}$} and convergence rate (c.r.) of CNFP method in solving wave problem (\ref{example1}) at time $t = 1$, where  $h = 1/64$,  and $s_1(x) = 1+0.3\sin(\pi x/8)$, and $s_2(x) = 1+0.2\tanh(\cos(\pi x/8))$.}
\label{table-CN}
\vspace{-3mm}
\end{center}
\end{table}
The results in Table \ref{table-CN} verify that the CNFP scheme has the second-order of accuracy in time for both constant-order and variable-order cases.  
Numerical simulations show that it generally takes longer time to solve the variable-order fractional wave equations due to computing the variable-order fractional Laplacian. 
It shows that the CNFP remains stable even with a large time step. 
At each time step the nonlinear system is solved by iterations, but the CNFP method enables us to use large time step,  which could potentially save computational costs by reducing the number of time steps to simulate. 

\subsection{Leap-frog Fourier pseudospectral method }
\label{section3-2}

The CNFP method has spectral accuracy in space and second-order accuracy in time. 
It is unconditionally stable, which allows large time steps in simulations. 
However,  the CNFP method is implicit, and iterations are required to solve the resulting system at each time step, which could complicate implementation and increase computing time. 
To avoid this, we present an explicit leap-frog method for the nonlinear fractional wave equation (\ref{model}). 

The leap-frog method is one of the most popular temporal discretization methods in solving the second-order wave equations. 
Using the leap-frog method to the semi-discretization problem in (\ref{semi}), we then obtain the leap-frog Fourier pseudospectral (LFFP) scheme as:
\begin{subequations}\label{LFscheme}
\begin{align}
&\frac{u_{\bf j}^{n+1}-2u_{\bf j}^n+u_{\bf j}^{n-1}}{\tau^2}=-\kappa (-\Dt)_h^{s(\bx_{\bf j})} u_{\bf j}^{n} + f(u_{\bf j}^{n}),  \quad \ \mbox{for} \  n = 1, 2, \ldots, \\
&u_{\bf j}^0 = \phi(\bx_{\bf j}), \\
&u_{\bf j}^1=\phi(\bx_{\bf j})+\tau\,\psi(\bx_{\bf j}) + \frac{\tau^2}{2} \bigg(-\kappa (-\Dt)_h^{s(\bx_{\bf j})}\phi(\bx_{\bf j}) + f\big(\phi(\bx_{\bf j})\big)\bigg), 
\end{align}
\end{subequations}
for ${\bf j} \in S_x$, where the initial conditions at $n = 0, 1$ are discretized exactly as in (\ref{IC-dis}). 
The LFFP method (\ref{LFscheme}) has the second-order temporal accuracy and spectral-order spatial accuracy, similar to the CNFP scheme in (\ref{CNscheme})--(\ref{IC-dis}). 
However, the LFFP scheme is fully explicit in time, and at each time step the computational cost is ${\mathcal O}(M N\log N)$ with $N$ the total number of spatial points. 
 
In Table \ref{table-LF},  we study the accuracy of the LFFP method in solving the fractional wave equation in (\ref{example1}), where the mesh size $h = 1/64$ is fixed. 
\begin{table}[htb!]
\begin{center}
\begin{tabular}{|c|c|c|c|c|c|c|c|c|c|c|}
\hline
\multirow{2}{*}{$\tau$} &\multicolumn{2}{|c|}{$s(x) \equiv 0.5$} & \multicolumn{2}{|c|}{$s(x) \equiv 1$} &  \multicolumn{2}{|c|}{$s(x) \equiv 1.3$} &  \multicolumn{2}{|c|}{$s_1(x)$} &  \multicolumn{2}{|c|}{$s_2(x)$} \\
\cline{2-11}
& error & c.r. &  & c.r. & error & c.r. & error & c.r. & error & c.r. \\
\hline
$2^{-7}$ & 1.135e-6 & --  &  5.079e-6 & -- & unstable & --  & unstable & -- & unstable & -- \\
$2^{-8}$ & 2.837e-7 & 2.00 &  1.269e-6 & 2.00 & unstable & --  & unstable & --  & 1.833e-6 & --  \\
$2^{-9}$ & 7.079e-8 & 2.00 &  3.167e-7 & 2.00 & 8.792e-7 & -- & 3.531e-7 & -- & 4.573e-7 & 2.00 \\
$2^{-10}$ & 1.756e-8 & 2.01 &  7.854e-8 & 2.01 & 2.181e-7 & 2.01 & 8.758e-8 & 2.01 & 1.134e-7 & 2.01 \\
$2^{-11}$ & 4.249e-9 & 2.04 &  1.901e-8 & 2.04 & 5.279e-8 & 2.04 & 2.120e-8 & 2.04 & 2.745e-8 & 2.04 \\
\hline
\end{tabular}
\caption{Temporal errors { $\|u(t) - u^{h, \tau}(t)\|_{l^2}$} and convergence rate (c.r.) of LFFP method in solving wave problem (\ref{example1}) at time $t = 1$, where  $h = 1/64$,  and $s_1(x) = 1+0.3\sin(\pi x/8)$, and $s_2(x) = 1+0.2\tanh(\cos(\pi x/8))$.}
\label{table-LF}
\vspace{-3mm}
\end{center}
\end{table}
The numerical parameters and reference solutions are prepared in the same manner as in Table \ref{table-CN}.  
It is clear that the LFFP scheme is conditionally stable. 
Table \ref{table-LF} shows that when the stability condition is satisfied,  it has the second-order temporal accuracy for both constant-order and variable-order cases.  
Compared to the implicit CNFP scheme, the explicit LFFP method is computationally much cheaper. 
Its computational efficiency becomes more significant in high dimensions. 

In the constant-order (i.e. $s(\bx) \equiv s$) cases, we can obtain the CFL condition for stability as $\tau < Ch^s$ with $C$ a positive constant independent of $h$. 
Particularly if $s = 1$, it collapses to the CFL condition of the leap-frog method in solving classical wave equations. 
It indicates that the smaller the value of $s$, the larger the threshold of time step, as usually $h < 1$.  
In the variable-order cases,  it is challenging to obtain the analytical stability condition as its constant-order counterpart.  
Denote $\tau^*$ as the critical time step for stability, i.e. the largest time step ensuring the method's stability.
\begin{figure}[htb!]
\centerline{
\includegraphics[width=6.86cm, height=5.06cm]{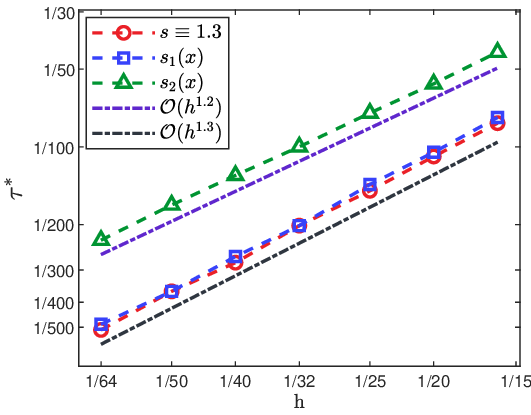}}
\caption{Critical time step $\tau^*$ versus mesh size $h$ for LFFP method, where $s_1(x) = 1+0.3\sin(\pi x/8)$, and $s_2(x) = 1+0.2\tanh(\cos(\pi x/8))$.}\label{Figure-2LFFP}
\end{figure}
Figure \ref{Figure-2LFFP} shows the relation between critical time step $\tau^*$ and  mesh size $h$, which verifies the CFL condition for constant $s$.
Moreover,  Figure \ref{Figure-2LFFP} suggests that in the variable-order cases, the critical time step $\tau^*$ depends on the maximum value of $s(x)$. 
Hence,  we introduce an enhanced CFL condition: 
\begin{equation}\label{enhancedCFL}
\displaystyle\tau <  Ch^{\max\{s(\bx)\}}. 
\end{equation}
Since the spectral method is applied in space, the mesh size $h>0$ usually can be rather large. 
This CFL condition is often acceptable in practice.   
Our extensive studies further confirm the stability condition (\ref{enhancedCFL}) for the LFFP method. 

\subsection{Time-splitting Fourier pseudospectral method}
\label{section3-3}

In this section, we propose another new explicit temporal discretization method. 
Denote $v(\bx,t):=\partial_t u(\bx,t)$, for $t \ge 0$. 
The second-order fractional wave equation in (\ref{model}) can be formulated into a first-order system of $(u, v)$, i.e.,
\begin{equation}\label{eq-order1}
\left\{\begin{split}
 &\partial_tu(\bx,t)=v(\bx,t), \\
 &\partial_tv(\bx,t)=-\kappa(-\Dt)^{s(\bx)}u(\bx,t)+f(u), \qquad \mbox{for} \ \ t > 0, 
 \end{split}\right.
\end{equation}
subject to the initial conditions 
$$u(\bx, 0) = \phi(\bx), \qquad  v(\bx, 0) = \psi(\bx). $$ 
At time $t = t_n$, assume solutions $u(\bx, t_n)$ and $v(\bx, t_n)$ are known.  
Then from  $t=t_n$ to $t=t_{n+1}$, we propose to split \eqref{eq-order1} and solve it in two steps, i.e., solving  
\begin{equation}\label{phik}
\left\{\begin{split}
 &\partial_t u(\bx,t)=v(\bx,t),\\
 &\partial_t v(\bx,t)=-\kappa(-\Dt)^{s_0}u(\bx,t),
 \end{split} \right. \quad \ \mbox{for} \ \ t \in(t_n,  t_{n+1}],  \hspace{3.4cm}
 \end{equation}
and
\begin{equation}\label{phip}
\left\{\begin{split}
 &\partial_t u(\bx,t)=0,\\
 &\partial_t v(\bx,t)=\kappa\left[(-\Dt)^{s_0}
 -(-\Dt)^{s(\bx)}\right]u(\bx,t) + f(u), 
 \end{split}\right.\quad \ \mbox{for} \ \ t \in(t_n,  t_{n+1}]. 
\end{equation} 

In the special case of constant order (i.e. $s(\bx) \equiv s$), we have $s_0 = s$. 
Consequently, the subproblem (\ref{phip}) reduces to a time-dependent ODE system, while the subproblem (\ref{phik}) is equivalent to the linear fractional wave equation (\ref{model}) with $s(\bx) \equiv s$ and $f(u) = 0$.

Next we focus on solving  (\ref{phik}) and (\ref{phip}).  
To this end, we assume that the solution $(u, v)$ at point $\bx = \bx_{\bf j}$ can be approximated by
\bea \label{phik-sol}
u_{\bf j}(t) = \sum_{{\bf k}\,\in\,{S}_\xi^0} \widehat{u}_{\bf k}(t)\,\fe^{i {\boldsymbol\mu_{\bf k} }{\boldsymbol\cdot}(\bx_{\bf j}-{\bf a})}, \quad \ 
v_{\bf j}(t) = \sum_{{\bf k}\,\in\,S_\xi^0} \widehat{v}_{\bf k}(t)\,\fe^{i {\boldsymbol\mu_{\bf k} }{\boldsymbol\cdot}(\bx_{\bf j}-{\bf a})}, \qquad\mbox{for ${\bf j} \in S_x$},
\eea
where $\widehat{u}_{\bf k}(t)$ and $\widehat{v}_{\bf k}(t)$ are defined in the same manner as in (\ref{ansatz-d}). 
Substituting (\ref{phik-sol}) into  (\ref{phik}) and taking FFT at both sides lead  to the following system: 
\beas
\fl{\md \widehat{u}_{\bf k}(t)}{\md t}  = \widehat{v}_{\bf k}(t), \quad \; \fl{\md \widehat{v}_{\bf k}(t)}{\md t}  = -\kappa\,|{\boldsymbol\mu}_{\bf k}|^{s_0}\,\widehat{u}_{\bf k}(t),  
 \qquad \mbox{for} \ \ {\bf k}\in {S}_\xi^0,
\eeas
which can be exactly integrated in time. 
We then obtain 
\begin{equation}\label{uhat1}
\begin{split}
&\widehat{u}_{\bf k}(t) = \widehat{u}_{\bf k}(t_n) \cos\big[w_{\bf k} (t-t_n)\big] + \fl{\widehat{v}_{\bf k}(t_n)}{w_{\bf k}}\sin\big[w_{\bf k}(t-t_n)\big], \\
&\widehat{v}_{\bf k}(t) = 
-w_{\bf k} \widehat{u}_{\bf k}(t_n) \sin\big[w_{\bf k} (t-t_n)\big]\,\widehat{u}_{\bf k}(t_n) + \widehat{v}_{\bf k}(t_n)\cos\big[w_{\bf k}(t-t_n)\big], 
\end{split} \quad \mbox{for}\ \  |{\bf k}| \neq  0, 
\end{equation}
 for $t \ge t_n$, where $w_{\bf k} = \sqrt{\kappa}\,|{\boldsymbol \mu}_{\bf k}|^{s_0}$. 
Note that if $|{\bf k}| = 0$, we have $|{\boldsymbol\mu}_{\bf k}| = 0$. 
Thus, we get 
\begin{equation}\label{uhat0}
\begin{split}
&\widehat{u}_{\bf k}(t) = \widehat{u}_{\bf k}(t_n) + (t-t_n)\,\widehat{v}_{\bf k}(t_n),\\
&\widehat{v}_{\bf k}(t) = \widehat{v}_{\bf k}(t_n),
\end{split} \qquad  \mbox{for}\ \  |{\bf k}| =  0,\qquad\qquad
\end{equation}
for $t \ge t_n$.
Combining (\ref{phik-sol})--(\ref{uhat0})  immediately gives the numerical solution of (\ref{phik}).
  
On the other hand, the sub-problem (\ref{phip}) can be integrated in time exactly and gives
\bea\label{phip-sol}
\begin{split}
&u_{\bf j}(t) = u_{\bf j}(t_n), \qquad \\
&v_{\bf j}(t) = v_{\bf j}(t_n) + (t-t_n) \Big(\kappa\widetilde{(-\Dt)}_h^{s(\bx_{\bf j})}u_{\bf j}(t_n)+ f\big(u_{\bf j}(t_n)\big)\Big),
\end{split}
\eea
for time $t \ge t_n$. 
Note that the term with $m = 0$ in (\ref{Lu-d})  coincides with the approximation of $(-\Dt)_h^{s_0}u_{\bf j}(t)$. 
Hence, we introduce the notation
\beas
\begin{aligned}
\widetilde{(-\Dt)}_{h}^{s(\bx_{\bf j})}u_{\bf j}(t) &:= (-\Dt)_h^{s(\bx_{\bf j})}u_{\bf j}(t)-\kappa (-\Dt)_h^{s_0}u_{\bf j}(t)\\
& \,\, = \sum_{m=1}^{M}\big[s({\bf x_j})- s_0\big]^m\Bigg(\sum_{{\bf k}\in{S_\xi}}\frac{\big(\ln{|\boldsymbol\mu_{\bf {\bf k}}|^2}\big)^m}{m!}\big[|\boldsymbol\mu_{\bf k}|^{2s_0}\widehat{u}_{\bf k}(t)\big]\fe^{i{\boldsymbol\mu}_{\bf k}{\boldsymbol\cdot}(\bx_{\bf j}-{\bf a})}\Bigg).
\end{aligned}
\eeas

From $t = t_n$ to $t = t_{n+1}$, we can combine  \eqref{phik} and \eqref{phip} by the second-order Strang splitting method and obtain the second-order time-splitting Fourier pseudospectral (TSFP2) method as follows:
\bea\label{splitscheme}
\begin{split}
&u_{\bf j}^{(1)} = \sum_{{\bf k}\,\in\,{S}_\xi^0} \widehat{u}_{\bf k}\big(t+\fl{\tau}{2}\big)\,\fe^{i {\boldsymbol\mu_{\bf k} }{\boldsymbol\cdot}(\bx_{\bf j}-a)}, \quad \;\;
v_{\bf j}^{(1)} = \sum_{{\bf k}\,\in\,{S}_\xi^0} \widehat{v}_{\bf k}\big(t+\fl{\tau}{2}\big)\,\fe^{i {\boldsymbol\mu_{\bf k} }{\boldsymbol\cdot}(\bx_{\bf j}-a)}, \\
&u_{\bf j}^{(2)} = u_{\bf j}^{(1)},  \quad \;\;
v_{\bf j}^{(2)} = v_{\bf j}^{(1)}+\tau\Big(\kappa\widetilde{(-\Dt)}_h^{s(\bx_{\bf j})}u_{\bf j}^{(1)}+ f\big(u_{\bf j}^{(1)}\big)\Big),\\
&u_{\bf j}^{n+1} = \sum_{{\bf k}\,\in\,{S}_\xi^0} \widehat{u}_{\bf k}^{(2)}\big(t+\fl{\tau}{2}\big)\,\fe^{i {\boldsymbol\mu_{\bf k} }{\boldsymbol\cdot}(\bx_{\bf j}-a)}, \quad \;\;
v_{\bf j}^{n+1} = \sum_{{\bf k}\,\in\,{S}_\xi^0} \widehat{v}_{\bf k}^{(2)}\big(t+\fl{\tau}{2}\big)\,\fe^{i {\boldsymbol\mu_{\bf k} }{\boldsymbol\cdot}(\bx_{\bf j}-a)}, \\
\end{split}\qquad 
\eea
for ${\bf j} \in S_{x}$ and $n = 1, 2, \ldots$,\, where $\widehat{u}_{\bf k}(t)$, $\widehat{v}_{\bf k}(t)$, $\widehat{u}_{\bf k}^{(2)}(t)$, and $\widehat{v}_{\bf k}^{(2)}(t)$ can be computed from (\ref{uhat1}) or (\ref{uhat0}). 
For $n = 0$, the initial conditions are {\it exactly} given by 
\bea\label{split-IC}
u_{\bf j}^0 = \phi(\bx_{\bf j}), \quad \; \ v_{\bf j}^0 = \psi(\bx_{\bf j}), \qquad \mbox{for} \ \ {\bf j} \in S_x.
\eea
The TSFP2 method in (\ref{splitscheme})--(\ref{split-IC}) has the second-order temporal accuracy and spectral-order spatial accuracy. 
It is explicit in time and easy to implement.  
The computational cost at each time step is ${\mathcal O}(MN\log N)$ with $N$ the total number of spatial points. 
Since two subproblems are integrated exactly in time, the only temporal errors in scheme (\ref{splitscheme}) are splitting errors, which can be improved by using higher order split step method \cite{Bao2005}. 

\begin{figure}[htb!]
\centerline{
\includegraphics[width=6.86cm, height=5.06cm]{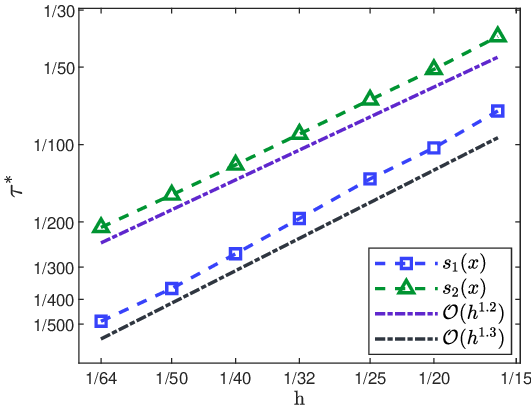}}
\caption{Critical time step $\tau^*$ versus mesh size $h$ for TSFP2 method, where $s_1(x) = 1+0.3\sin(\pi x/8)$, and $s_2(x) = 1+0.2\tanh(\cos(\pi x/8))$.}\label{Figure2-TSFP2}
\end{figure}

In Table \ref{table-TS},  we numerically study the accuracy of the TSFP2 method in solving the fractional wave equation (\ref{example1}), where the numerical parameters and reference solutions are prepared in the same manner as those in Tables \ref{table-CN}--\ref{table-LF}. 
\begin{table}[htb!]
\begin{center}
\begin{tabular}{|c|c|c|c|c|c|c|c|c|c|c|}
\hline
\multirow{2}{*}{$\tau$} &\multicolumn{2}{|c|}{$s(x) \equiv 0.5$} & \multicolumn{2}{|c|}{$s(x) \equiv 1$} &  \multicolumn{2}{|c|}{$s(x) \equiv 1.3$} &  \multicolumn{2}{|c|}{$s_1(x)$} &  \multicolumn{2}{|c|}{$s_2(x)$} \\
\cline{2-11}
& error & c.r. &  & c.r. & error & c.r. & error & c.r. & error & c.r. \\
\hline
$2^{-7}$ & 7.661e-6 & --  &  5.893e-6 & -- & 5.096e-6 & --  & unstable & -- & unstable & -- \\
\hline
$2^{-8}$ & 1.914e-6 & 2.00 &  1.473e-6 & 2.00 & 1.274e-6 & 2.00 & unstable& -- & 1.020e-6 & -- \\
\hline
$2^{-9}$ & 4.776e-7 & 2.00 &  3.674e-7 & 2.00 & 3.180e-7 & 2.00 & 3.876e-7 & -- & 2.544e-7 & 2.00 \\
\hline
$2^{-10}$ & 1.185e-7 & 2.01 &  9.113e-8 & 2.01 & 7.904e-8 & 2.01 & 9.612e-8 & 2.01 & 6.310e-8 & 2.01 \\
\hline
$2^{-11}$ & 2.868e-8 & 2.04 &  2.206e-8 & 2.04 & 1.931e-8 & 2.03 & 2.327e-8 & 2.04 & 1.527e-8 & 2.04 \\
\hline
\end{tabular}
\caption{Temporal errors $\|u(t) - u^{h, \tau}(t)\|_{l^2}$ and convergence rate (c.r.) of TSFP2 method in solving  wave problem (\ref{example1}) at time $t = 1$, where  $h = 1/64$,  and $s_1(x) = 1+0.3\sin(\pi x/8)$, and $s_2(x) = 1+0.2\tanh(\cos(\pi x/8))$.}
\label{table-TS}
\vspace{-3mm}
\end{center}
\end{table}
The results in Table \ref{table-TS} confirm the second-order temporal accuracy of the TSFP2 method.  
In contrast to the LFFP method, the TSFP2 method is stable in solving constant-order fractional wave equations. 
While solving variable-order wave equations, TSFP2 method has a similar stability condition as LFFP method; see Figure \ref{Figure2-TSFP2}. 
Hence,   the stability condition in (\ref{enhancedCFL}) also applies to the TSFP2 method. 
The TSFP2 method is explicit and thus especially efficient in solving high-dimensional wave equations. 
Moreover, the temporal accuracy can be further increased by using the high-order Strang splitting method \cite{Bao2005, Duo2016}.

\begin{remark}
In the special case of constant-order linear fractional wave equation, the problem reduces only to solve sub-problem (\ref{phik})  and no temporal errors are introduced.  
\end{remark}

\section{Numerical experiments}
\label{section4}

In this section, we present numerical experiments, on one hand,  to further study the performance of our methods, and on the other hand to study the properties of fractional wave equations. 
To the best of our knowledge, numerical report on the fractional wave equations with variable-order fractional Laplacian still remains very limited due to the lack of effective methods for approximating the variable-order fractional Laplacian. 
In existing studies \cite{Zhu2014, Yao2017}, the constant-order fractional Laplacians are used to approximate the variable-order Laplacian, which significantly reduces the computational complexity.  
However, this essentially changes heterogeneous media to homogeneous media and could strongly affect the wave propagation and dynamics. 

Unless otherwise stated, we will use the proposed matrix-free approach (\ref{Lu-d}) to calculate the variable-order fractional Laplacian in the following studies.

\subsection{Accuracy verification}
\label{section4-1}

We test the numerical accuracy of three numerical methods, including CNFP in (\ref{CNscheme})--(\ref{IC-dis}), LFFP in (\ref{LFscheme}), and TSFP2 in (\ref{splitscheme}). 
To this end, we solve the following one-dimensional nonlinear fractional wave equation: 
\bea\label{example1}
\begin{split}
& \p_{tt}u(x, t) = -\kappa(-\p_{xx})^{s(x)}u + u^3, \quad \ \  \mbox{for} \ \ t > 0, \\
&u(x, 0) = \exp(-x^2), \qquad \p_tu(x, 0) = 0.
\end{split}
\eea
In our simulations,  we choose $\kappa = 1$ and set the computational domain as $[-32, 32]$.  
The exact solution of this nonlinear wave problem is unknown.  
When computing numerical errors,  we use the numerical solution with fine mesh size $h = 1/64$ and small time step $\tau = 0.0001$ as the reference ``exact" solution. 
Let $u(t)$ denote this reference ``exact" solution at time $t$, while $u^{h, \tau}(t)$ represents the numerical solution computed with mesh size $h$ and time step $\tau$. 
In the following, we choose  $M = 15$ which is large enough such that the truncation errors of matrix-free approach do not affect the spatial and temporal errors of our methods. 

In Tables \ref{table-CN}--\ref{table-TS},  the temporal accuracy of CNFP, LFFP, and TSFP2 methods have been studied, where we fix the mesh size $h = 1/64$. 
It shows that all three methods have the second-order of accuracy in time. 
The CNFP method is implicit, and at each time step iterations are required to solve the resultant system. 
Hence, the computing time for CNFP is much longer than that for LFFP and TSFP methods; see detailed comparison in Table \ref{Table-CPU}. 
Next, we test the spatial accuracy.  
In Table \ref{Table-spatialerror}, we only show the spatial errors of the  TSFP2 method since the same spatial discretization is used for all three methods. 
Here,  the small time step $\tau = 0.0001$ is fixed.
\begin{table}[htb!]
\begin{center}
\begin{tabular}{|c|c|c|c|c|c|}
\hline
$h$&{$s(x) \equiv 0.5$} & {$s(x) \equiv 1$} &  {$s(x) \equiv 1.3$} &  {$s_1(x)$} & {$s_2(x)$} \\
\hline
$1$ & 2.7614e-2 &  4.7112e-2 &  3.6899e-2 &  5.1481e-2 &  2.8845e-2  \\
\hline
$1/2$ &  6.4113e-4 &  2.3151e-4 &  2.4355e-4 &  5.4924e-4 &  1.9511e-4  \\
\hline
$1/4$ &  6.5436e-7 &  2.8530e-8 &  1.6516e-9 &  2.7637e-7 &  5.8982e-9  \\
\hline
$1/8$ &  2.449e-12 &  5.460e-13 &  6.014e-13 &  5.958e-13 &  5.637e-13  \\
\hline
\end{tabular}
\caption{Spatial discretization errors $\|u(t) - u^{h, \tau}(t)\|_{l^2}$ of TSFP method at time $t = 1$, where time step $\tau = 0.0001$. 
The variable order $s_1(x) = 1+0.3\sin(\pi x/8)$, and $s_2(x) = 1+0.2 {\rm tanh}\big[\cos(\pi x/8)\big]$. 
Note that similar spatial errors are obtained for CNFP and LFFP methods.}
\label{Table-spatialerror}
\end{center}
\vspace{-3mm}
\end{table}
It is evident that the spatial discretization has a spectral order accuracy. 
Our extensive studies show that CNFP and LFFP methods have the similar spatial errors which we will omit showing for brevity. 

\begin{figure}[htb!]
\centerline{
\includegraphics[height = 5.06cm, width=6.86cm]{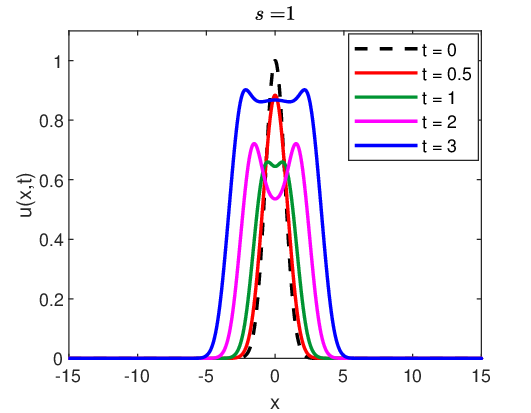}\hspace{-4mm}
\includegraphics[height = 5.06cm, width=6.86cm]{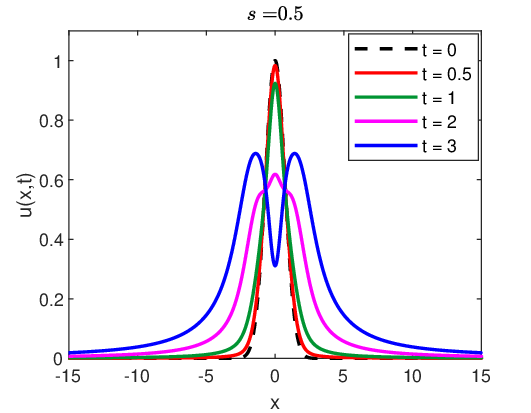}}
\centerline{
\includegraphics[height = 5.06cm, width=6.86cm]{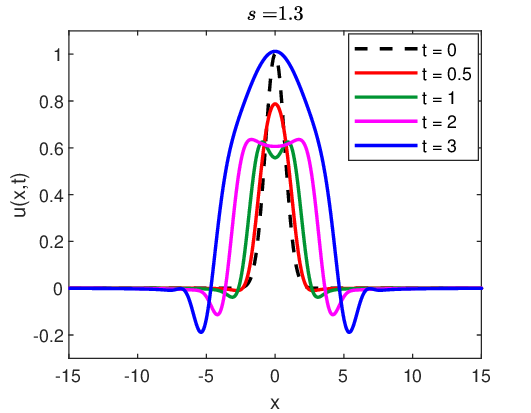}\hspace{-4mm}
\includegraphics[height = 5.06cm, width=6.86cm]{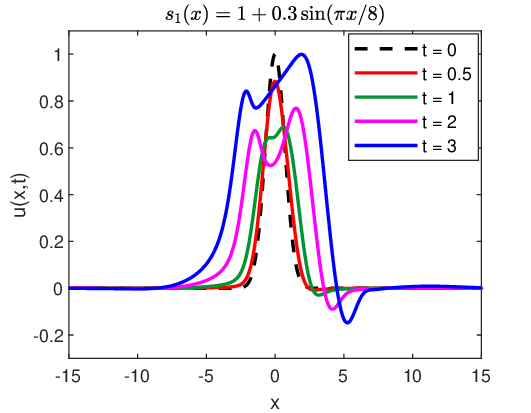}}
\caption{Solution dynamics of the nonlinear  wave equation in (\ref{example1}) with different $s(x)$. 
For better illustration, the displayed domain $[-15, 15]$ is much smaller than our computational domain.}\label{Figure-Ex1}
\end{figure}

Figure \ref{Figure-Ex1} presents the solution dynamics for different $s(x)$. 
It shows that in all cases the solution disperses,  and grows over time due to the nonlinear term.  
Compared to classical ($s = 1$) cases, the nonlocality of the fractional Laplacian leads to a faster dispersion. 
In homogeneous (i.e. constant $s$) media, the solution remains symmetric for any time $t \ge 0$. 
In contrast, it becomes asymmetric in heterogeneous media. 
The solution dynamics in heterogeneous media are more complicated due to the interplay of nonlocality and heterogeneity. 
Importantly, we note that the dynamics of waves for $s_1(x)$ is very different from that in homogeneous media with $s(x) \equiv 1$, even though the average value of $s_1(x)$ is equal to 1. 
This suggests that using averaged constant-order fractional Laplacians to approximate the variable-order Laplacian fails to describe the heterogeneity of wave propagation. 

\subsection{Efficient comparison} 
\label{section4-2}

We continue to test and compare the performance of CNFP, LFFP, and TSFP2 methods, especially focusing on the effectiveness of our matrix-free approach  (\ref{Lu-d}) for spatial approximation.  
To this end, we focus on the spatially varying $s(x)$. 
The nonlinear fractional wave equation in  (\ref{example1}) is solved, where computational domain $\Og = [-32, 32]$ and time step $\tau = 0.0001$. 

\begin{figure}[ht!]
\centerline{
(a)\includegraphics[height=5.56cm, width=6.8cm]{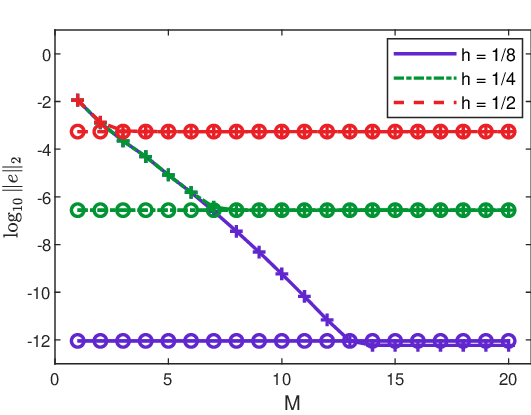}\hspace{2mm}
(b)\includegraphics[height=5.46cm, width=6.8cm]{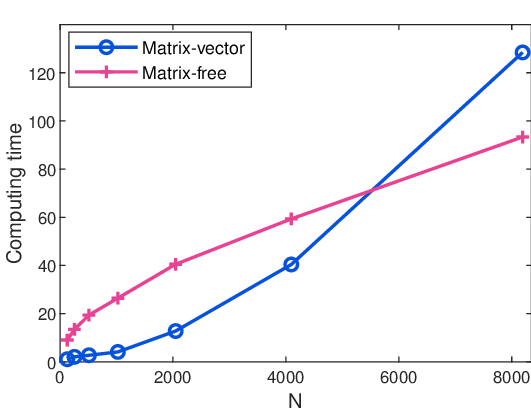}
}
\caption{Comparison of the matrix-free approach (symbol `$+$') and the direct matrix-vector approach (symbol `$\circ$') in solving the problem (\ref{example1}) with $s(x) = 1+0.3\sin(\pi x/8)$, where the TSFP2 method is used with $\tau = 0.0001$. (a) Numerical errors at time $t = 1$; (b) Computing time taken from $t = 0$ to $t = 1$.}\label{Figure-Error}
\end{figure}

First, we study the truncation effects of $M$ in matrix-free approach by comparing it to the direct matrix-vector approach in (\ref{FL1Ddis0}). 
Figure \ref{Figure-Error} (a) compares the numerical errors of these two approaches for different $h$ and $M$. 
Here, the time-splitting method is adopted for temporal discretization, and the reference ``exact" solution is prepared in the same way as in Example 1. 
As expected, numerical errors of the matrix-vector approach (\ref{FL1Ddis0}) are independent of $M$. 
Note that our matrix-free approach can be viewed as a truncated approximation of the direct matrix-vector method, and its truncation error depends on  $M$.   
Figure \ref{Figure-Error} (a) shows that numerical errors of the matrix-free approach decrease exponentially as $M$ increases, eventually converging to those of the matrix-vector approach. 
This is because the truncation errors of matrix-free approach are dominant when $M$ is small.
However, when $M$ is large enough,  the discretization errors of numerical methods become prominent, and the truncation errors can be neglected. 
Moreover, the threshold value of $M$ depends on the mesh size $h$. 
These observations suggest that the matrix-free scheme with a large $M$ can accurately compute the variable-order fractional Laplacian while significantly reducing the computational and storage burden caused by its heterogeneity. 

In Figure \ref{Figure-Error} (b), we compare the computing time of these two approaches for solving the one-dimensional problem (\ref{example1}), where we take $M = 15$ in the matrix-free approach (\ref{Lu}). 
It shows that when the number of points $N$ is small, the matrix-free scheme takes longer. 
However, the computing time of the direct matrix-vector method drastically increases as $N$ grows larger. 
At each time step, the computational cost of direct method is ${\mathcal O}(N^2)$, whereas the matrix-free approach costs ${\mathcal O}(M N\log N)$. 
Moreover, the storage cost of the direct approach is ${\mathcal O}(N^2)$, in contrast to ${\mathcal O}(MN)$ of the matrix-free method. 
The memory limitation could become the main bottleneck of the direct matrix-vector approach in high dimensions. 
To see it, Table \ref{Tab-2D} presents their computing time in solving two-dimensional fractional wave equations, where $M = 20$ in the matrix-free approach. 
\begin{table}[htb!]
   \centering
    \begin{tabular}{|c|c|c|c|c|c|c|} \hline
      $N$     &$16^2$  &$32^2$   &$64^2$    &$128^2$  &$256^2$ &$512^2$  \\ \hline
    matrix-vector scheme       &0.925   &3.023   &27.58   &379.1 & n.a. & n.a.\\ \hline
    matrix-free scheme       &19.57    &39.14   &79.48   &191.8 &560.3 &2202 \\ \hline
    \end{tabular}
    \caption{Computing time of the matrix-vector and  matrix-free approaches in solving the two-dimensional wave equation with $s(\bx)=1-0.4\cos({\pi x}/{4})\cos({\pi y}/{4})$ and $t\in(0,1]$, where the TSFP2 method is used with $\tau = 0.0001$.} \label{Tab-2D}
\end{table}
Our studies are conducted on a laptop equipped with an Intel(R) Core(TM) i7-12700H processor and 32GB of RAM.
It shows that in two-dimensional cases, storing the entire matrix for $N \ge 256^2$ becomes impossible, making the matrix-vector method infeasible. 
Therefore, the advantages of the matrix-free scheme are more pronounced in higher dimensions. 
Furthermore, the outer summation in (\ref{Lu-d}) can be easily parallelized, further reducing computing time and enhancing efficiency

Our extensive studies show that the same conclusions drawn for TSFP2 method from Figures \ref{Figure-Error} and Table \ref{Tab-2D} can be applied to CNFP and LFFP methods.  
In Table \ref{Table-CPU}, we further compare the computing time of CNFP, LFFP, and TSFP2 methods in solving the one-dimensional fractional wave equation for $t \in (0, 1]$ with time step $\tau = 0.0001$. 
\begin{table}[ht]
\centering
\begin{tabular}{|c|r|r|r|r|r|r|r|r|} \hline
\multirow{2}{*}{$N$}&\multicolumn{2}{|c|}{CNFP method} &\multicolumn{2}{|c|}{LFFP method} &\multicolumn{2}{|c|}{TSFP2 method}   \\ \cline{2-7}
              & m.-v.& m.-f.  & m.-v.&  m.-f.  & m.-v.& m.-f.\\ \hline
            1024 &20.672 &100.667   &1.097   &25.539 & 4.070  & 26.555 \\ 
            2048 &290.537 &{155.591}  &3.945  &{ 39.710} & 11.406 & 41.224 \\ 
            4096 &1113.809 &240.702 &28.021 &57.497 & 45.799 & 59.933 \\ 
            8192 &4850.179 &503.139 &129.998 &87.996 & 132.960& { 93.341} \\ \hline
\end{tabular}
\caption{Computing time of the matrix-vector (m.-v.) and  matrix-free (m.-f.) approaches in solving 1D nonlinear fractional wave equation (\ref{example1}) for $t\in(0, 1]$, where $s(x)= 1+0.3\sin(\pi x/8)$ and time step $\tau = 0.0001$.}
        \label{Table-CPU}
\end{table}
Consistent to the observations in Figure \ref{Figure-Error} (b), the accelerated matrix-free approach takes more time than the direct matrix-vector method for small number of points $N$. 
Compared to the explicit methods, CNFP takes much longer time as at each time step the nonlinear system is iteratively solved, which becomes more problematic if the matrix-vector scheme is used.

\subsection{Application simulations}
\label{section4-3}

In the following, we numerically study the solution dynamics of fractional wave equations and compare the nonlocal effects of homogeneous and heterogeneous fractional Laplacians. 

\bigskip
\noindent{\bf Example 1 (Soliton collision). } 
We numerically study the interaction of two { solitary waves}  so as to understand the nonlocality and/or heterogeneity of the (variable-order) fractional Laplacian. 
For this purpose,  the one-dimensional linear wave equations  is considered with the initial conditions
\bea\label{example3}
\begin{split}
    u(x, 0)&= {\rm sech}{\big(a(x+x_0)\big)} +{\rm sech}{\big(a(x-x_0)\big)}, \\
    u_t(x,0) &= \fl{b\sinh{\big(a(x+x_0)\big)}}{\cosh^2{\big(a(x+x_0)\big)}}-\fl{b\sinh{\big(a(x-x_0)\big)}}{\cosh^2{\big( a(x-x_0)\big)}},
\end{split}
\eea
where we choose $a=b=3$ and $x_0 = 10$.  
Initially, two well-separated waves are centered at $x = \pm 10$, respectively.  
The coefficient $\kappa$ is set as $\kappa = 1$. 
Choose the computational domain $[-128, 128]$ with mesh size $h = 1/64$,  and time step $\tau = 0.0001$. 
We have verified that our results are independent of numerical parameters by refining both mesh size $h$ and time step $\tau$.

\begin{figure}[ht]
\centerline{
\includegraphics[width=5.16cm, height=4.46cm]{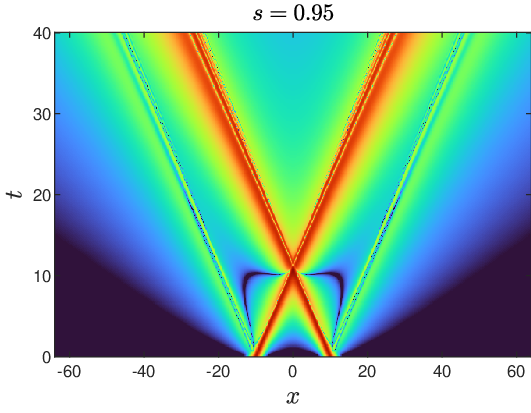} 
\includegraphics[width=5.16cm, height=4.46cm]{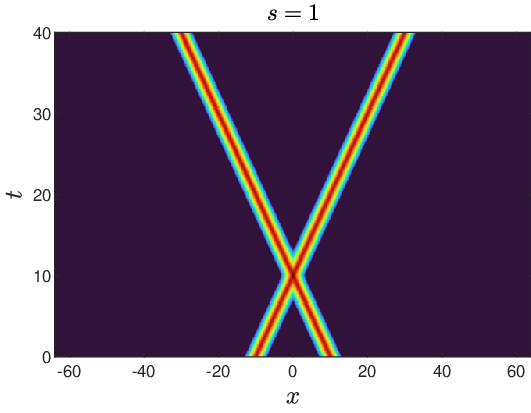}
\includegraphics[width=5.16cm, height=4.46cm]{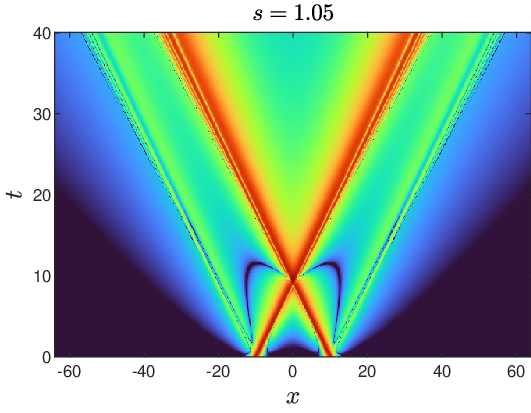}}
\centerline{
\includegraphics[width=5.26cm, height=4.26cm]{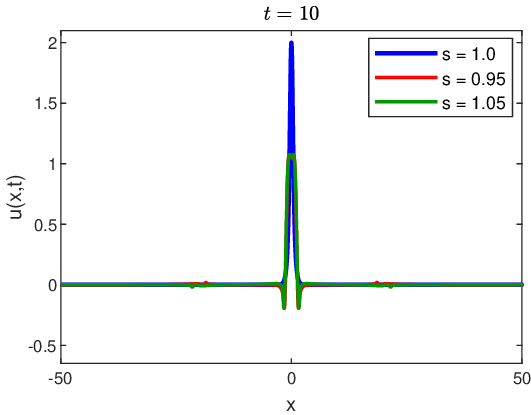}\hspace{-1mm}
\includegraphics[width=5.26cm, height=4.26cm]{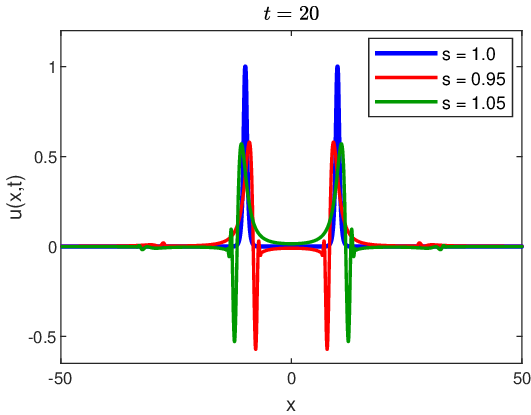}\hspace{-1mm}
\includegraphics[width=5.26cm, height=4.26cm]{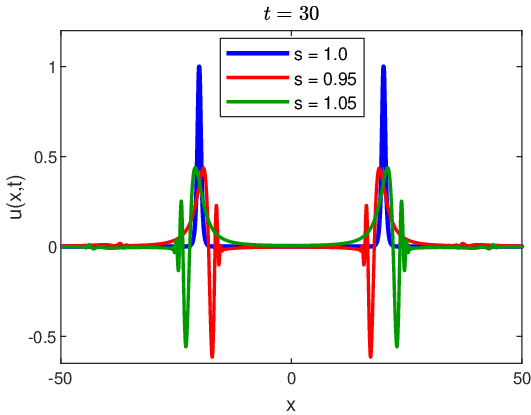}}
\caption{Dynamics of two solitons in homogeneous media with $\kappa = 1$ and $f = 0$ in (\ref{example1}), and initial condition in (\ref{example3}). For better illustration, the displayed domain is much smaller than our computational domain. }\label{Figure-soliton1}
\end{figure}

In Figure \ref{Figure-soliton1}, we illustrate the dynamics of two solitons in homogeneous media (i.e., constant $s$).  
It shows that two initially well-separated solitons first move towards each other, and then collide at $x = 0$ around time $t = 10$. 
After collision, they separate again and move apart. 
In the classical ($s \equiv 1$) case,   two solitons can retain their initial shape after separation. 
However,  two solitons in the fractional cases continue to deform over time,  and moreover radiation of waves is observed. 
Even though the two waves move at an equal speed in homogeneous media, their speed depends on the value of $s$ -- the larger the value of $s$, the faster the wave moves (cf. Figure \ref{Figure-soliton1} for $s = 0.95$ and $1.05$). 
Figure \ref{Figure-soliton1} further shows that  the radiation patterns for $s < 1$ and $s > 1$ are different.  
Our extensive studies show that when $s$ is further away from $1$, the solution dynamics become more chaotic, consistent with the observations in \cite{Kirkpatrick2016}.

\begin{figure}[ht!]
\centerline{
a)\includegraphics[width=6.06cm, height=4.6cm]{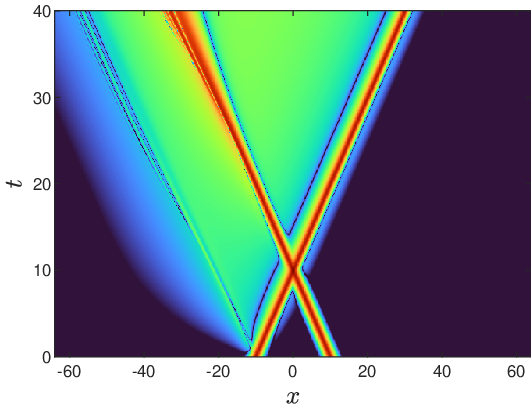} \hspace{3mm}
b)\includegraphics[width=6.06cm, height=4.6cm]{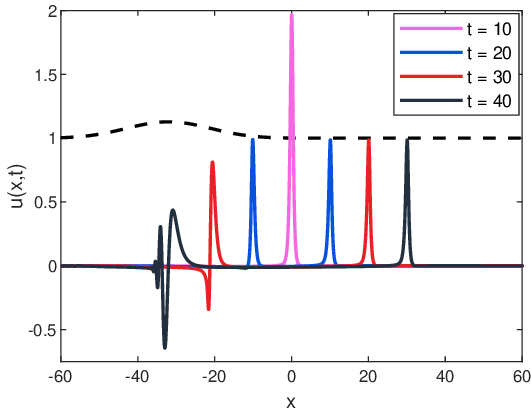}}\vspace{5mm}
\centerline{
c)\includegraphics[width=6.06cm, height=4.6cm]{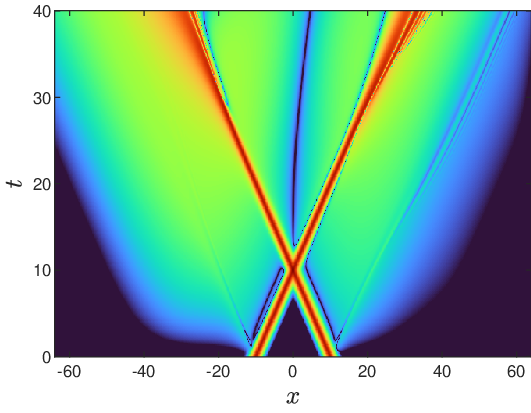} \hspace{3mm}
d)\includegraphics[width=6.06cm, height=4.6cm]{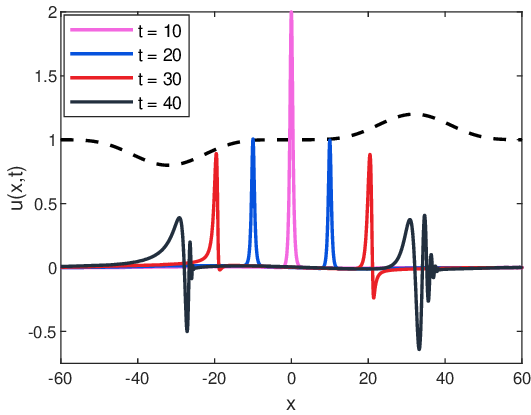}}
\caption{Dynamics of two solitons in heterogeneous media with $\kappa = 1$ and $f = 0$ in (\ref{example1}), and initial condition in (\ref{example3}), where the dashed line in b) \& d) represents $s(x)$. 
a)--b): $s(x) = 1+0.001\big[1-\sin(\pi x/64)\big]^7$; c)--d): $s(x) = 1+0.2\sin^5(\pi x/64)$.  For better illustration, the displayed domain is much smaller than our computational domain.  }\label{Figure-soliton2}
\end{figure}

By comparison,  Figure \ref{Figure-soliton2} shows the dynamics of two waves in heterogeneous media with $s_3(x) = 1+0.001\big[1-\sin(\pi x/64)\big]^7$ and $s_4(x) = 1+0.2\sin^5(\pi x/64)$. 

It shows that two waves collide and subsequently move apart, similar to those observed in Figure \ref{Figure-soliton1}. 
However, the dynamics of two waves in this case are asymmetric. 
Let's use $s_3(x)$ as an example, which describes a heterogeneous medium with classical medium for $x \ge 0$ and fractional medium for $x < 0$. 
From Figure \ref{Figure-soliton2} a) and b),  we find that the two waves collide around $t = 10$, and then they separate and move apart.
After separation, the wave on the right-hand side maintains its original profile as moving in the classical ($s = 1$) media. 
While the wave on the left-hand side exhibits qualitative behavior as that observed in Figure \ref{Figure-soliton1} c) for $s = 1.05$. 
It is clear that the evolution of two waves is asymmetric due to the heterogeneity characterized by the variable-order fractional Laplacian. 
The case of $s_4(x)$ can be viewed as a composite comprising three media. 
Waves' radiation is observed on both sides over a long time. 

\bb
\noindent{\bf Example 2 (Wave propagation). } 
We  study the wave dispersion in the two-dimensional fractional wave equation to further understand the nonlocal effect of the (variable-order) fractional Laplacian.  
Consider  the two-dimensional ($d = 2$) linear fractional wave equation (\ref{model}) with $\kappa = 0.2$. 
The initial conditions are taken as 
\beas
u(\bx, 0) = 5\Big(e^{-20(x^2+(y+0.1)^2)}- e^{-20(x^2+(y-0.1)^2)}\Big), \qquad u_t(\bx, 0) = 0.
\eeas
In our simulations, we take the computational domain as $[-12, 12]^2$ with number of points $J_1 = J_2 = 2048$.  
The time step is $\tau = 0.0001$. 

\begin{figure}[ht!]
\centerline{
\includegraphics[width=3.46cm, height=3.06cm]{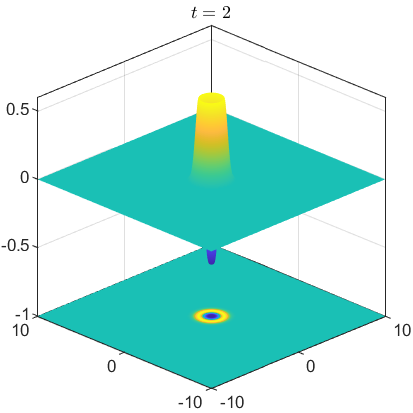}
\includegraphics[width=3.46cm, height=3.06cm]{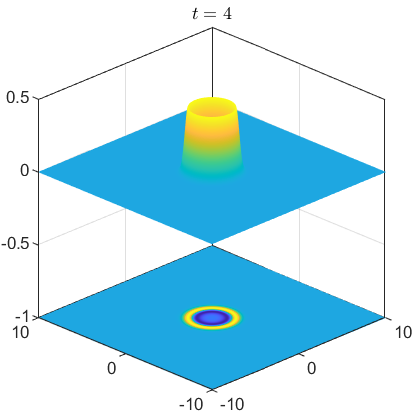}
\includegraphics[width=3.46cm, height=3.06cm]{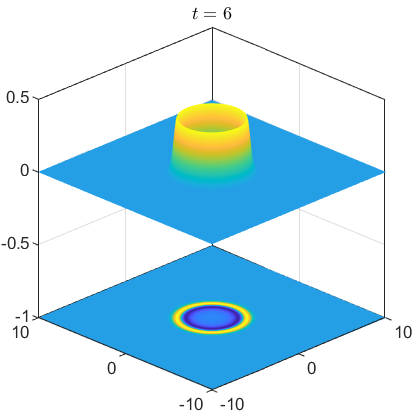}
\includegraphics[width=3.46cm, height=3.06cm]{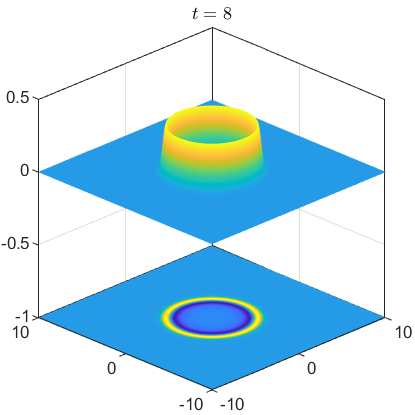}}
\vspace{2mm}
\centerline{
\includegraphics[width=3.46cm, height=3.06cm]{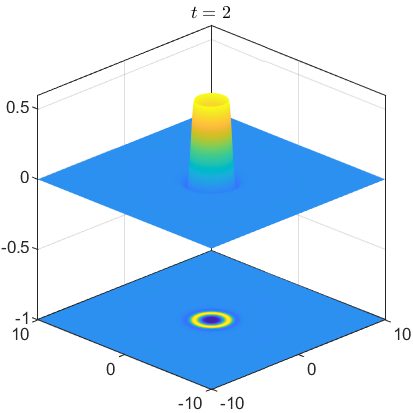}
\includegraphics[width=3.46cm, height=3.06cm]{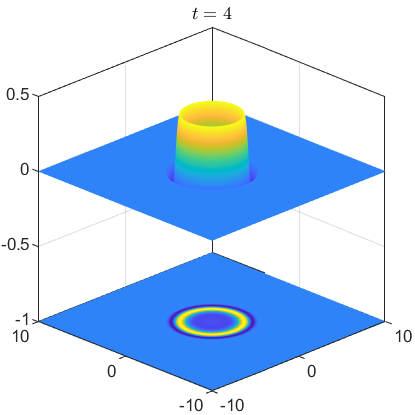}
\includegraphics[width=3.46cm, height=3.06cm]{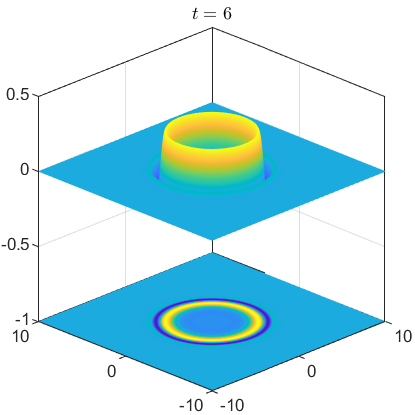}
\includegraphics[width=3.46cm, height=3.06cm]{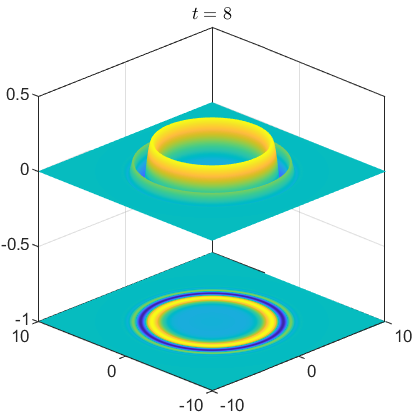}}
\vspace{2mm}
\centerline{
\includegraphics[width=3.46cm, height=3.06cm]{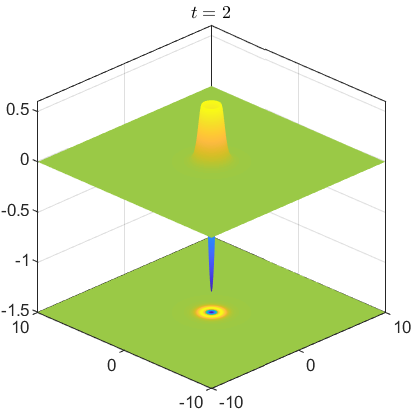}
\includegraphics[width=3.46cm, height=3.06cm]{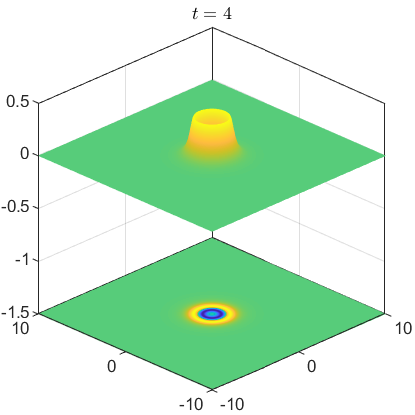}
\includegraphics[width=3.46cm, height=3.06cm]{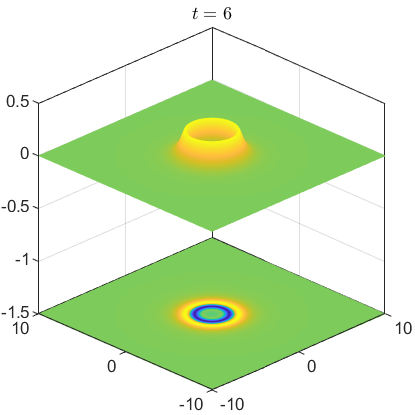}
\includegraphics[width=3.46cm, height=3.06cm]{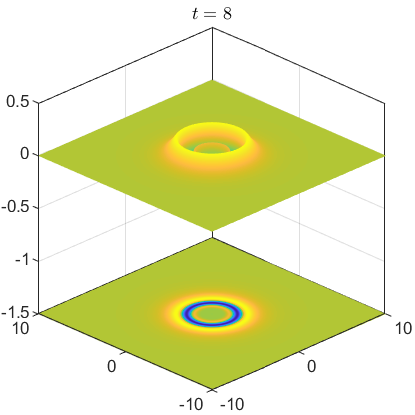}}
\caption{Solution dynamics in two-dimensional wave equations in homogeneous media, where $s = 1, 1.2$, and $0.8$ from top to bottom rows.}\label{Figure-2Dconstant}
\end{figure}

\begin{figure}[ht!]
\centerline{
\includegraphics[width=3.46cm, height=3.06cm]{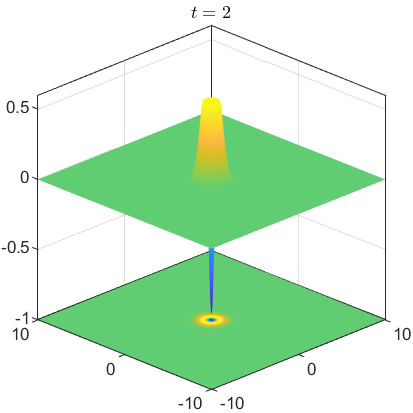}
\includegraphics[width=3.46cm, height=3.06cm]{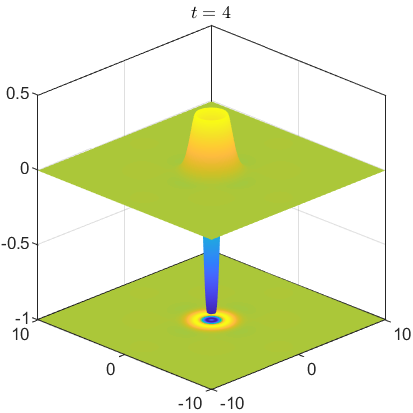}
\includegraphics[width=3.46cm, height=3.06cm]{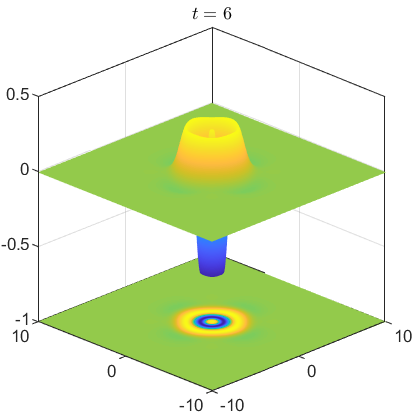}
\includegraphics[width=3.46cm, height=3.06cm]{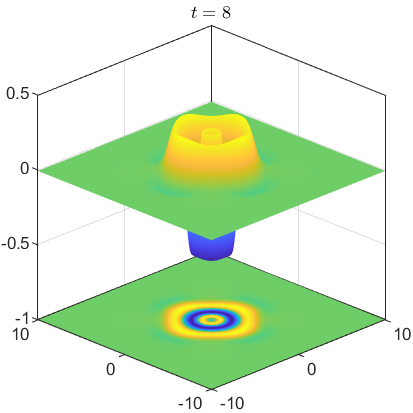}}
\vspace{2mm}
\centerline{
\includegraphics[width=3.46cm, height=3.06cm]{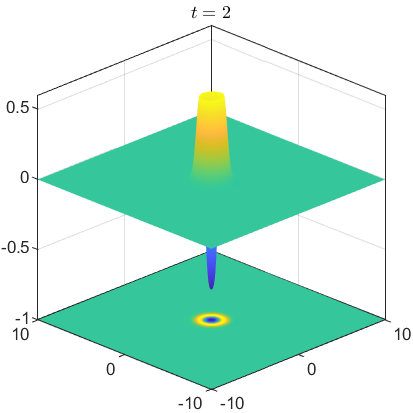}
\includegraphics[width=3.46cm, height=3.06cm]{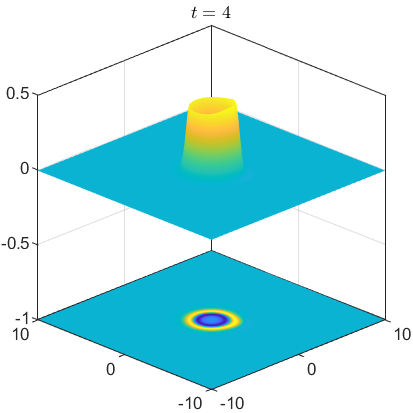}
\includegraphics[width=3.46cm, height=3.06cm]{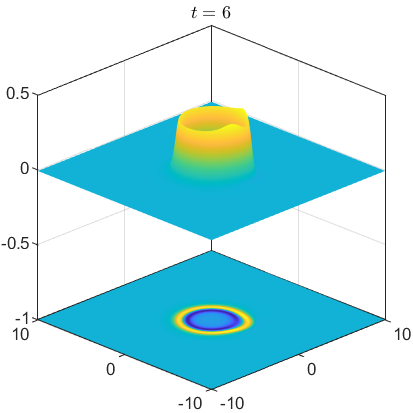}
\includegraphics[width=3.46cm, height=3.06cm]{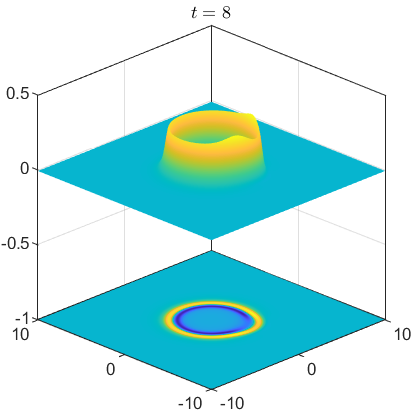}}
\vspace{2mm}
\caption{Solution dynamics in two-dimensional wave equations in heterogeneous media, where $s(\bx) = 1-0.4\cos(\fl{\pi x}{4})\cos(\fl{\pi y}{4})$ (top) and $s(\bx) = 1+0.2\exp{(-(x-1)^2-(y+2)^2)}$ (bottom).}\label{Figure-2Dvariable}
\end{figure}

Figure \ref{Figure-2Dconstant} shows the solution dynamics in constant-order fractional wave equations, where the results of the classical ($s \equiv 1$) wave equations are included as a benchmark. 
It shows that the solution of constant $s$ propagates radially outward over time. 
The larger the constant $s$, the faster the expansion of solution. 
If $s \neq 1$, the solution scatters over time. 
The scattering occurs towards the center for $s < 1$, and towards the boundary for $s > 1$. 
By contrast, Figure \ref{Figure-2Dvariable} illustrates the solution dynamics of variable-order fractional wave equation. 
Because of the heterogeneity, the solutions evolve asymmetrically, and their velocities in different directions depend on the medium parameter $s(\bx)$. 
Consistent with the observations in Figure \ref{Figure-2Dconstant}, solution propagates more rapidly in the region with large value of $s(\bx)$. 
Moreover, our results support the conclusion in \cite{Meerschaert2014} that the variable-order models enable one to capture the anisotropic nature of wave propagation in complex media. 
Computationally, it is more challenging to study the fractional wave equation in heterogeneous media due to the spatially varying $s(\bx)$. 
Our numerical studies show that our matrix-free approach is effective, and the outer summation in (\ref{Lu-d})  can be parallelized to further reduce computing time. 

\bb
\noindent{\bf Example 3 (Wave dispersion and attenuation). } 
We extend our method to study wave dispersion and attenuation in a simple two-layer heterogeneous attenuating media \cite{Zhu2014, Yao2017}. 
To this end, we consider the 
two-dimensional fractional wave equation in (\ref{seismicwave}). 
Choose  $\gamma(\bx) = a_1 + a_2{\rm tanh}\big(100(y-1)\big)$, i.e. representing a two-layered media separated at $y = 1$. 
The initial conditions take the form 
\begin{equation}\label{seismic-IC}
    \phi(\bx) = \big( 1-2\pi^2 \nu_{0}^2|\bx-\bx_{c}|^2\big){\rm e}^{-\pi^2\nu_{0}^2|\bx-\bx_{c}|^2}, \qquad \psi(\bx,0)=0,
\end{equation}
where $\nu_0 = 25$ and $\bx_c = (1, 0.85)^T$ are chosen in our study.

Set the computational domain $\Og = (0, 2)^2$. 
\begin{figure}[ht!]
\centerline{
\includegraphics[width=3.1cm, height=3.16cm]{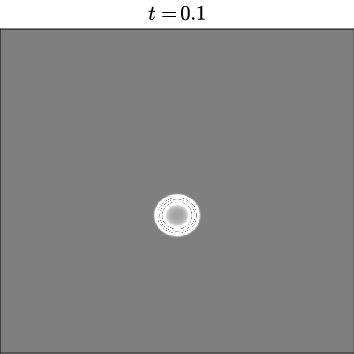}
\includegraphics[width=3.1cm, height=3.16cm]{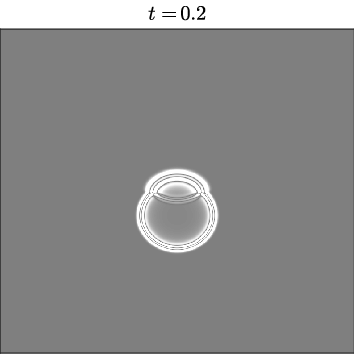}
\includegraphics[width=3.1cm, height=3.16cm]{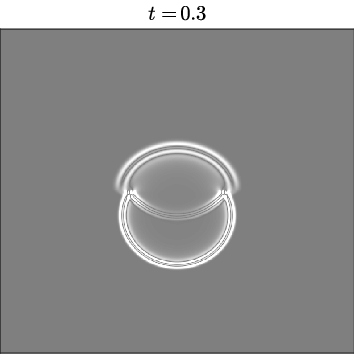}
\includegraphics[width=3.1cm, height=3.16cm]{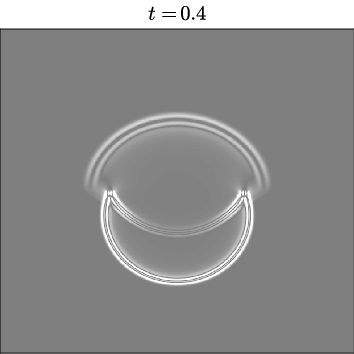}
\includegraphics[width=3.1cm, height=3.16cm]{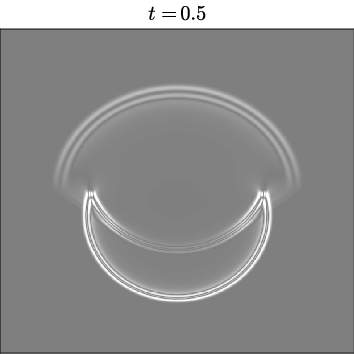}}
\centerline{
\includegraphics[width=3.1cm, height=3.16cm]{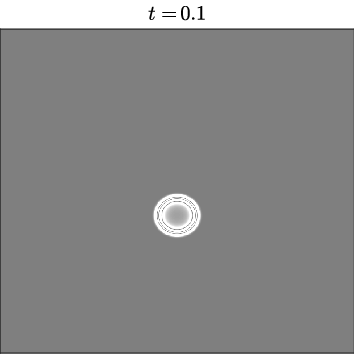}
\includegraphics[width=3.1cm, height=3.16cm]{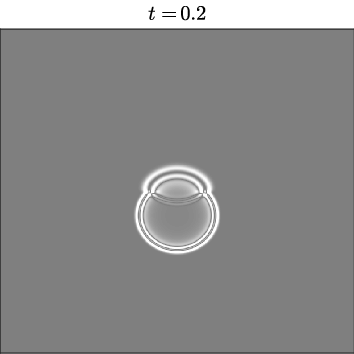}
\includegraphics[width=3.1cm, height=3.16cm]{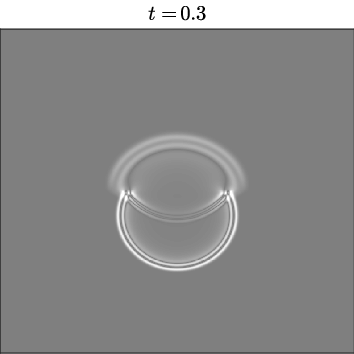}
\includegraphics[width=3.1cm, height=3.16cm]{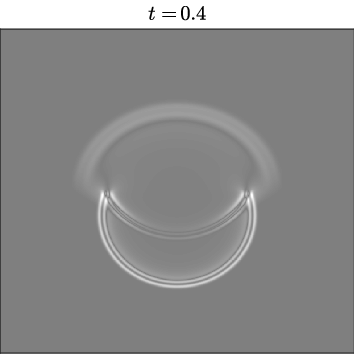}
\includegraphics[width=3.1cm, height=3.16cm]{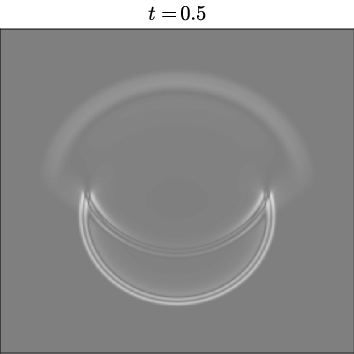}}
\caption{Wave propagation in two-layer attenuating media, where  $a_1 = 0.0065$ and $a_2 = 0.0035$ (first row), and $a_1 = 0.02$ and $a_2 = 0.01$ (second row).}\label{seismic1-variable}
\end{figure}
The model (\ref{seismicwave})  is solved by the LFFP method with mesh size $h = 0.002$ and time step $\tau = 0.0001$. 
Figure \ref{seismic1-variable} depicts wave propagation in heterogeneous media with different $\gamma(\bx)$, where $c_0 = 11/18$ and $1$ for the upper ($y > 1$) and lower ($y \le 1$) layers, respectively.
The dispersion and attenuation of waves are observed over time.  
The strength of attenuation increases with larger values of $\gamma(\bx)$, aligning with findings in \cite{Zhu2014, Yao2017}. 
Moreover, the wave decoherence is observed during the dynamics.

\section{Conclusions}
\label{section5}

We proposed accelerated Fourier pseudospectral methods to solve the variable-order space fractional wave equation. 
The spatial discretization is realized by the Fourier pseudospectral method, while the Crank--Nicolson, leap-frog, and time splitting methods are introduced for temporal discretization. 
In the special case of constant order (i.e. $s(\bx) \equiv s$),  our methods can be efficiently implemented via the (inverse) fast Fourier transforms, and the computational cost at each time step is ${\mathcal O}(N\log N)$, where $N$ represents the total number of spatial points.  
However, this fast algorithm fails in the variable-order cases due to the spatial dependence of the Fourier multiplier. 
To address this, we proposed an accelerated matrix-free approach for the efficient computation of variable-order cases. 
The computational cost is ${\mathcal O}(MN\log N)$ and storage cost ${\mathcal O}(MN)$, where $M \ll N$.  
Moreover, our scheme can be easily parallelized to further enhance its efficiency.  
It can be also applied to solve other variable-order fractional problems.

Numerical experiments were reported to examine the effectiveness of our methods and study the wave dynamics in heterogeneous media. 
Our numerical studies showed that in high dimensions,  the direct matrix-vector multiplication approach becomes impractical due to excessive memory requirements. 
In contrast, our accelerated method has proven effective, leveraging modest storage and computational resources. 
It showed that all three (CNFP, LFFP, and TSFP2) methods have the second-order accuracy in time. 
The explicit LFFP and TSFP2  methods require less computing time at each time step, while the CNFP method allows larger time step. 
We found that the wave dynamics in the fractional cases are more complicated due to the nonlocal interactions, especially in heterogeneous media compared to homogeneous media. 
The solution behaviors of the fractional wave equation will be further explored in our future study.

\section*{Acknowledgements}
X. Zhao is partially supported by the Natural Science Foundation of Hubei Province No. 2019CFA007 and the NSFC 11901440. 
Y. Zhang is partially supported by the US National Science Foundation DMS--1913293 and DMS--1953177.

\end{document}